\documentstyle[11pt,leqno,here]{article}
\input amssym.def
\input amssym.tex
\input psfig

\newcommand{\beql}[1]{\begin{equation}\label{#1}}
\newcommand{\eeq}{\end{equation}}
\newcommand{\EE}{{\Bbb E}}

\makeatletter
\@addtoreset{figure}{section}
\def\thefigure{\thesection.\@arabic\c@figure}
\makeatother

\begin{document}
\begin{center}
{\Large {\bf Dense Packings of Equal Disks in an Equilateral \\
\vspace{.1\baselineskip}
Triangle: From 22 to 34 and Beyond}} \\
\vspace{\baselineskip}
{\em R. L. Graham} \\
{\em B. D. Lubachevsky} \\
\vspace{.25\baselineskip}
AT\&T Bell Laboratories, \\
Murray Hill, New Jersey 07974 \\
\vspace{1\baselineskip}
\vspace{1.5\baselineskip}
{\bf ABSTRACT}
\vspace{.5\baselineskip}
\end{center}

\setlength{\baselineskip}{1.5\baselineskip}

Previously published packings of equal disks in an equilateral
triangle have dealt with up to 21 disks.
We use a new discrete-event simulation algorithm to produce
packings for up to 34 disks.
For each $n$ in the range $22 \le n \le 34$ we present what we believe
to be the densest possible packing of $n$ equal disks in an
equilateral triangle.
For these $n$ we also list the second, 
often the third and sometimes the fourth best packings 
among those that we found.
In each case, the structure of the packing implies that the minimum
distance $d(n)$ between disk centers is the root of polynomial
$P_n$ with integer coefficients.
In most cases we do not explicitly compute $P_n$ 
but in all cases we do compute and report $d(n)$ 
to 15 significant decimal digits.

Disk packings in equilateral triangles differ 
from those in squares or circles in that for triangles 
there are an infinite number of values of $n$ 
for which the exact value of $d(n)$ is known,
namely, when $n$ is of the form $\Delta (k) := \frac{k(k+1)}{2}$.
It has also been conjectured that $d(n-1) = d(n)$ in this case.
Based on our computations, we present conjectured optimal packings
for seven other infinite classes of $n$, namely
\begin{eqnarray*}
n & = & \Delta (2k) +1,~\Delta (2k+1) +1,
\Delta (k+2) -2 , ~ \Delta (2k+3) -3, ~ \\
&& \Delta (3k+1)+2 ,
~ 4 \Delta (k), ~~\mbox{and}~~
2 \Delta (k+1) + 2 \Delta (k) -1 ~.
\end{eqnarray*}
We also report the best packings we found 
for other values of $n$ in these forms 
which are larger than 34, namely, 
$n=37$, 40, 42, 43, 46, 49, 56, 57, 60, 63, 67, 71, 79, 84, 92, 93, 106, 112, 121, and 254,
and also for $n=58$, 95, 108, 175, 255, 256, 258, and 260.
We say that an infinite class of packings of $n$ disks,
$n=n(1), n(2),...n(k),...$, 
is {\em tight }, if
[$1/d(n(k)+1) - 1/d(n(k))$] is
bounded away from zero as $k$ goes to infinity.
We conjecture that some of our infinite classes are tight,
others are not tight, and that there are infinitely many 
tight classes.
\section{Introduction}
\hspace*{\parindent}
Geometrical packing problems have a long and distinguished history in combinatorial mathematics.
In particular, such problems are often surprisingly difficult.
In this note, we describe a series of computer experiments
designed to produce dense packings of equal nonoverlapping disks 
in an equilateral triangle.
It was first shown by Oler in 1961 \cite{O} 
that the densest packing of 
$n= \Delta (k) : = \frac{k(k+1)}{2}$ equal disks is the appropriate
triangular subset of the regular hexagonal packing of the
disks (well known to pool players in the case of $n=15$).
It has also been
conjectured by Newman \cite{N} (among others) 
that the optimal packing of $\Delta (k)-1$ disks 
is always obtained by removing a single disk 
from the best packing for $\Delta (k)$,
although this statement has not yet been proved.
The only other values of $n$ (not equal to $\Delta(k)$) 
for which optimal packings are known are
$n$ = 2, 4, 5, 7, 8, 9, 11 and 12 
(see Melissen \cite{M1}, \cite{M2} for a survey).

\begin{figure}[htb]
\centerline{\psfig{file=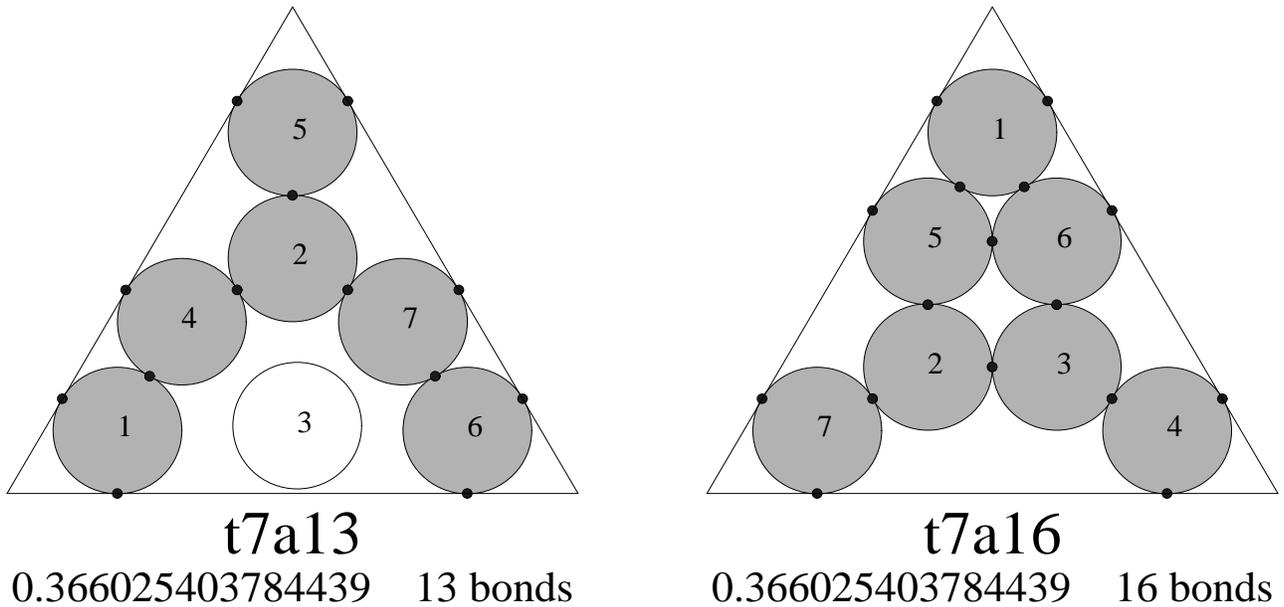,width=7in}}

\caption{Two equivalent but nonisomorphic densest packings of 7 disks.}
\label{t7}
\end{figure}

As the number $n$ of packed disks increases, it becomes
not only more difficult to {\em prove} optimality of a packing
but even to {\em conjecture} what the optimal packing might be.
In this paper, we present a number of conjectured optimal packings.
These packings are produced on a computer using a so-called
``billiards'' simulation algorithm.
A detailed description of the philosophy, 
implementation and applications of this event-driven
algorithm can be found in \cite{L}, \cite{LS}.
Essentially, the algorithm simulates a system of $n$ perfectly elastic disks.
In the absence of gravitation and friction, 
the disks move along straight lines, 
colliding with each other and the region walls 
according to the standard laws of mechanics, 
all the time maintaining a condition of no overlap.
To form a packing, 
the disks are uniformly allowed to gradually increase in size, 
until no significant growth can occur.
Not infrequently, it can happen at this point that
there are disks which can still move,
e.g., disk 3 in t7a13 (see Fig.~\ref{t7}).

Every packing of $n$ disks occurring
in the literature for $n$ different from $\Delta(k)$
and $\Delta(k)-1$
which has been conjectured or proved
to be optimal was also found by our algorithm.
These occur for
$n=13$, 16, 17, 18, and 19 (see \cite{M1}, \cite{MS}).
This increases our confidence that
the new packings we obtain are also optimal.
The new packings cover two ``triangular periods'':
$21 = \Delta (6)$ to $\Delta (7)$ to $\Delta (8) =36$.

In addition, we conjecture optimal packings for seven infinite classes of $n$,
namely, $n= \Delta (2k) +1$, 
$\Delta (2k+1) +1$, $\Delta (k+2) -2$, $\Delta (2k+3)-3$,
$\Delta (3k+1)+2$, $4 \Delta (k)$, and $2 \Delta (k+1) + 2 \Delta (k)-1$,
where $k=1,2...$.
Each class has its individual pattern of the optimal packings
which is different from patterns for other classes.
These were suggested by the preceding packings, 
and we give packings for some additional values of these forms, namely, 
$n=37$, 40, 42, 43, 46, 49, 56, 57, 60, 63, 67, 71, 79, 84, 92, 93, 106, 112, 121, and 254, as well as for $n=58$, 95, 108, 175, 255, 256, 258, and 260.

We say that an infinite class of packings of $n$ disks,
$n=n(1), n(2),...n(k),...$, 
is {\em tight }, if
[$1/d(n(k)+1) - 1/d(n(k))$] is
bounded away from zero as $k$ goes to infinity.
We conjecture that some of our infinite classes are tight,
others are not, and that there are infinitely many 
tight classes.
\section{The packings}
\hspace*{\parindent}
We performed a small number of runs with $n=21$, 27, 28, 35 and 36 disks.
In every case, the resulting packings were consistent with the
existing results ($n= \Delta (k)$) and conjectured
$(n= \Delta (k) -1)$.
The bulk of our efforts concentrated on the other 11 values of $n$,
for $21 \le n \le 36$.
These are presented in Figures~\ref{t22} to \ref{t34}.

To navigate among the various packings presented we will use the labeling
system illustrated by Fig.~\ref{t22}~t22a.
Here, $n=22$, ``a'' denotes that the packing is the best we found,
``b'' would be the second best (as in t23b in Fig.~\ref{t23}),
``c'' would be third best, and ``d'' would be fourth best.

Small black dots in the packing diagrams
are ``bonds''
whose number is also entered by each packing.
For example, there are 47 bonds in t22a.
A bond between two disks or between a disk
and a boundary indicates that the distance
between them is zero.
The {\em absence} of a bond in a spot where disk-disk or disk-wall
are apparently touching each other means that the corresponding
distance is strictly {\em positive},
though perhaps too small for the resolution of
the drawing to be visible.
For example, there is no bond between disk 1 and the left side
of the triangle in t18a (Fig.~\ref{t18});
according to our computations,
the distance between disk 1 and the side is 0.0048728...
of the disk diameter.
(Packing t18a was constructed in \cite{M1}.)
Each disk in most of the packings is provided with a label
which uniquely identifies the disk in the packing.
This labeling is nonessential;
it is assigned in order to facilitate referencing.
\begin{figure}[htb]
\centerline{\psfig{file=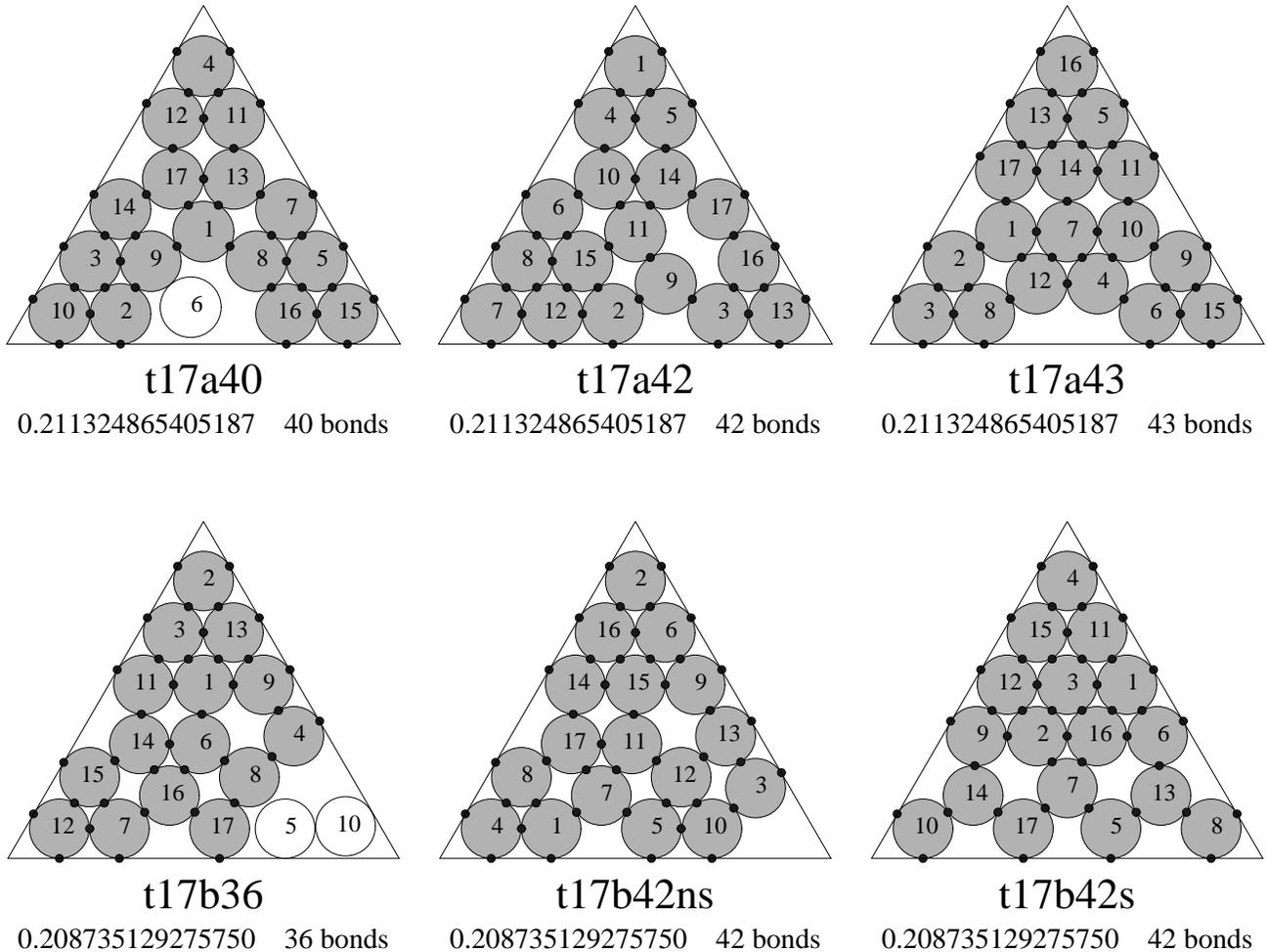,width=7in}}
\caption{The best (t17a40, t17a42, t17a43) and the next-best (t17b36, t17b42ns, t17b42s) packings of 17 disks.}
\label{t17}
\end{figure}
\begin{figure}[htb]
\centerline{\psfig{file=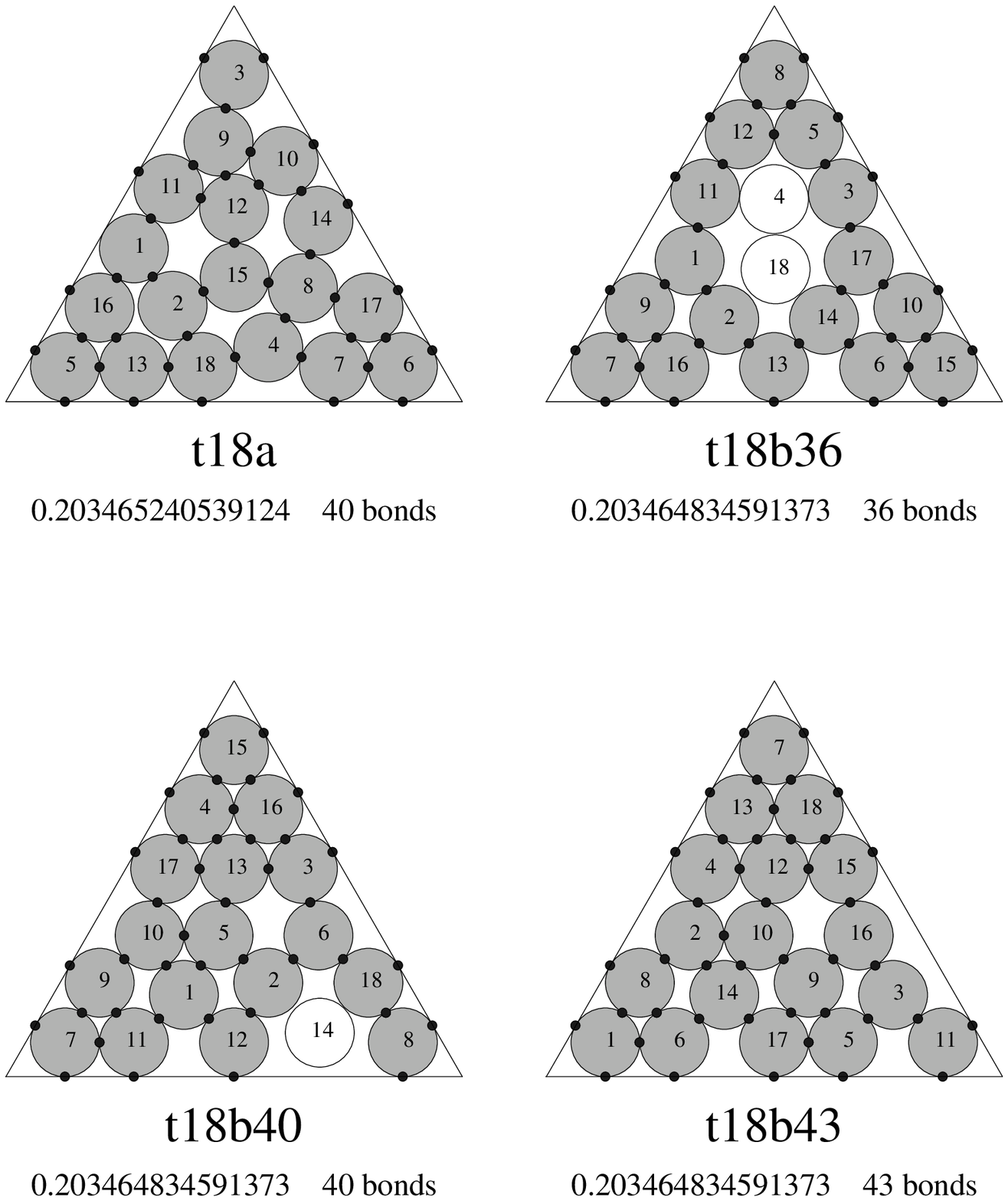,width=7in}}

\caption{The best (t18a) and the next best (t18b36, t18b40, t18b43) packings of 18 disks.}
\label{t18}
\end{figure}

Each disk normally has at least three bonds
attached.
The polygon formed by these bonds as vertices
contains the center of the disk
strictly inside.
This is a necessary condition for packing ``rigidity''.
In \cite{LS}, where the packing algorithm was applied
to a similar problem,
the disks without bonds were called ``rattlers.''
A rattler can move freely within the confines
of the ``cage'' formed by its
rigid neighbors and/or boundaries.
(If we ``shake'' the packing,
the rattler will ``rattle'' while hitting its cage.)
t22a has two rattlers, disks 3 and 5.
In the packing diagrams, all disks, except for the rattlers,
are shaded.

A number with 15 significant digits 
is indicated for each packing in the figures,
e.g., the number $0.17939~69086~11866$ for packing t22a.
This number is the disk diameter $d(n)$ which is measured in units
equal to the side of the smallest equilateral
triangle that contains the {\em centers} of all disks.
For packing t22a such a triangle 
is the one with vertices at the centers
of disks 22, 17, and 12.
This unit of measure for $d(n)$
conforms with previously published conventions.

Sometimes several packings exist
for the same disk diameter.
An example is t7a13 and t7a16 in Fig.\ref{t7}.
Thus, we
distinguish such packings by suffixing their labels
with the number of bonds.
Other examples are
t17a40, t17a42, and t17a43 in Fig.~\ref{t17},
t22b42 and t22b50 in Fig.~\ref{t22}.
However, even the number of bonds may not distinguish
different packings of the same disk diameter;
for example, t17b42ns and t17b42s in Fig.~\ref{t17},
where the provisional ``ns'' stands for ``non-symmetric''
and ``s'' for ``symmetric.''

We point out that the a-packings of 17 and 18 disks
that we show have previously been given by
Melissen and Schuur \cite{MS}, who also conjecture their optimality.
\section{Additional comments}
\vspace*{-.1in}
\paragraph{Fig.~\protect\ref{t23}:}
Two more c-packings for 23 disks 
that are not shown in the figure were generated:
t23c55.1 and t23c55.2. Both have 55 bonds.
t23c55.1  can be obtained by combining the left side
of t23c53 with the right side of t23c57.
t23c55.2 is a variant of t23c55.1.

\addtolength{\textheight}{.25in}
\paragraph{Fig.~\protect\ref{t24}:}
Disk 20 in t24c56 and in t24c59 is locked
in place because
its center is strictly inside
the triangle formed by the three bonds of disk 20.
In both packings, the distance of the disk center to the boundary 
of this enclosing triangle is the distance to the
line between bonds with the left side of
the triangle and disk 24, and is 0.0317185...
of the disk diameter.

\paragraph{Fig.~\protect\ref{t25}:}
The given d-packing of 25 disks t25d60 is symmetric
with respect to the vertical axis. An equivalent non-symmetric d-packing
t25d53 was also obtained in which all disks are located in the same
places as in t25d60, except for disks 5, 12, 13, 14, 23, and 24.
These six disks form a pattern which is roughly equivalent
to that formed by disks
10, 14, 19, 25, 20, and 22, respectively, in t25b.
Disk 24 in t25d53 is a rattler.

\paragraph{Fig.~\protect\ref{t29}:} 
Only one of the two b-packings of 29 disks we found is shown,
namely, t29b63.2.
The other b-packing, t29b63.1,
differs in the placements of only disks
2, 3, 4, 7 as
explained in Section~4.

\paragraph{Fig.~\protect\ref{t31}:}
Four a-packings of 31 disks exist;
only three are shown in the figure; the fourth one, t31a81.1,
is described in Section~4.

\paragraph{Fig.~\protect\ref{t33}:}
In t33a, the gap between disk 8 and left side
is 0.0017032... of the disk diameter.
In t33c, disk 7 is stably locked by its bonds
with 3, 6, and 29.
However, the distance from disk 7 center to the line on bonds
with disks 3 and 6 is only 0.0002097575.... of the disk diameter.
As a result, the cage of rattler disk 5 in t33c is very tight:
the gap between disk 22 and disk 5 or disk 18 and disk 5
does not exceed $4 \times 10^{-9}$ of the disk diameter.

\paragraph{Fig.~\protect\ref{t34}:}
In t34a, the small gaps between
``almost'' touching pairs disk-disk or disk-wall 
take on only three values
(relative to the disk size):
in pairs 20--31, 16--26, 23--27, 18--19, 1--27 the gap is 0.021359...,
in pairs left-32, right-29 it is 0.024750...,
and in pairs 4--34, 7--22, it is 0.042561...
Similarly, 
there are only three values of gaps in each of t34b, t34c, and t34d.
\begin{quote}
t34b: in pairs 18--19, 23--27, 17--28, 20--31, 16--26 the gap is 0.019583...;
in pairs left-32, right-29 it is 0.022686...; in pairs 4--34, 7--22, it is 0.039035... \\
t34c: in pairs 12--17, 22--27, 3--10, 14--21, 4--34. 3--16 the gap is 0.018864...;
in pair left-15 it is 0.021850...; in pair 19-24 it is 0.037606... \\
t34d: in pairs 2--4, 26--32, 15--22, 12--21, 3--16, 7--16 the gap is 0.018681...;
in pair left-27 it is 0.021637...; in pairs 13--33, 19--30 it is 0.037242... 
\end{quote}
\begin{figure}[htb]
\centerline{\psfig{file=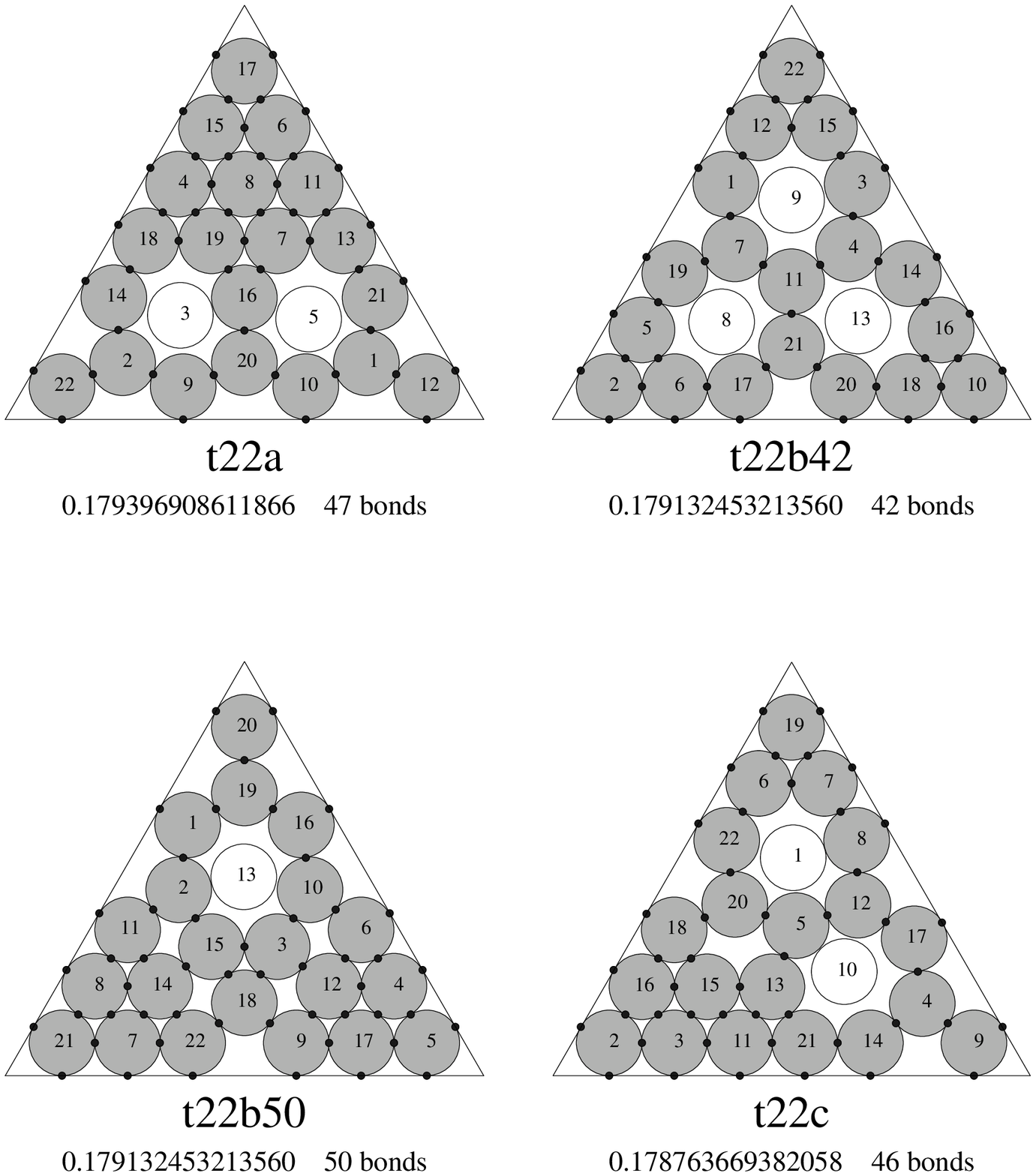,width=7in}}

\caption{The best (t22a), the next-best (t22b42, t22b50), and the third-best (t22c) packings of 22 disks.}
\label{t22}
\end{figure}
\addtolength{\textheight}{-.25in}
\begin{figure}[htb]
\centerline{\psfig{file=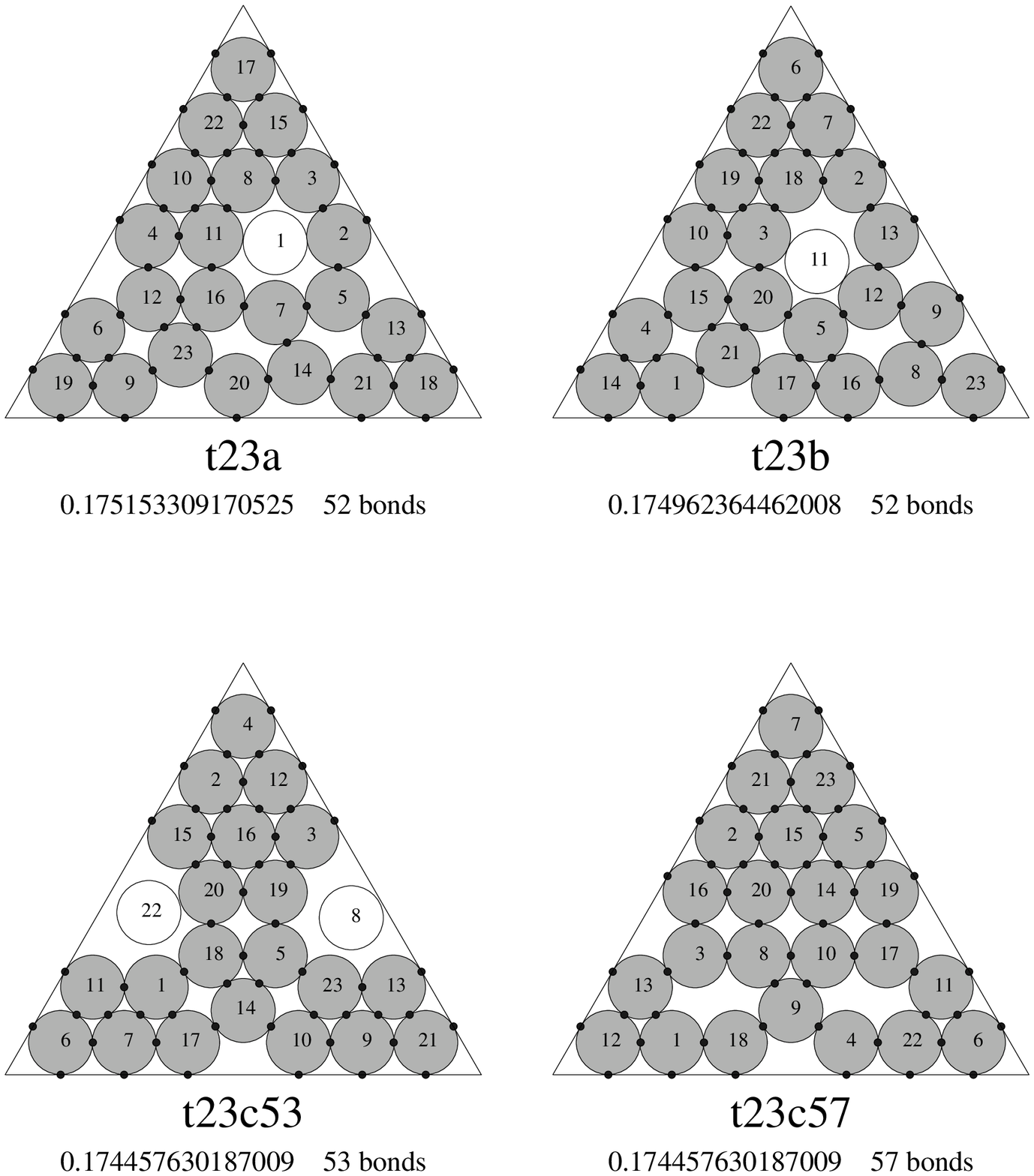,width=7in}}

\caption{The best (t23a), the next-best (t23b), and the third-best (t23c53, t23c57) packings of 23 disks.}
\label{t23}
\end{figure}
\begin{figure}[htb]
\centerline{\psfig{file=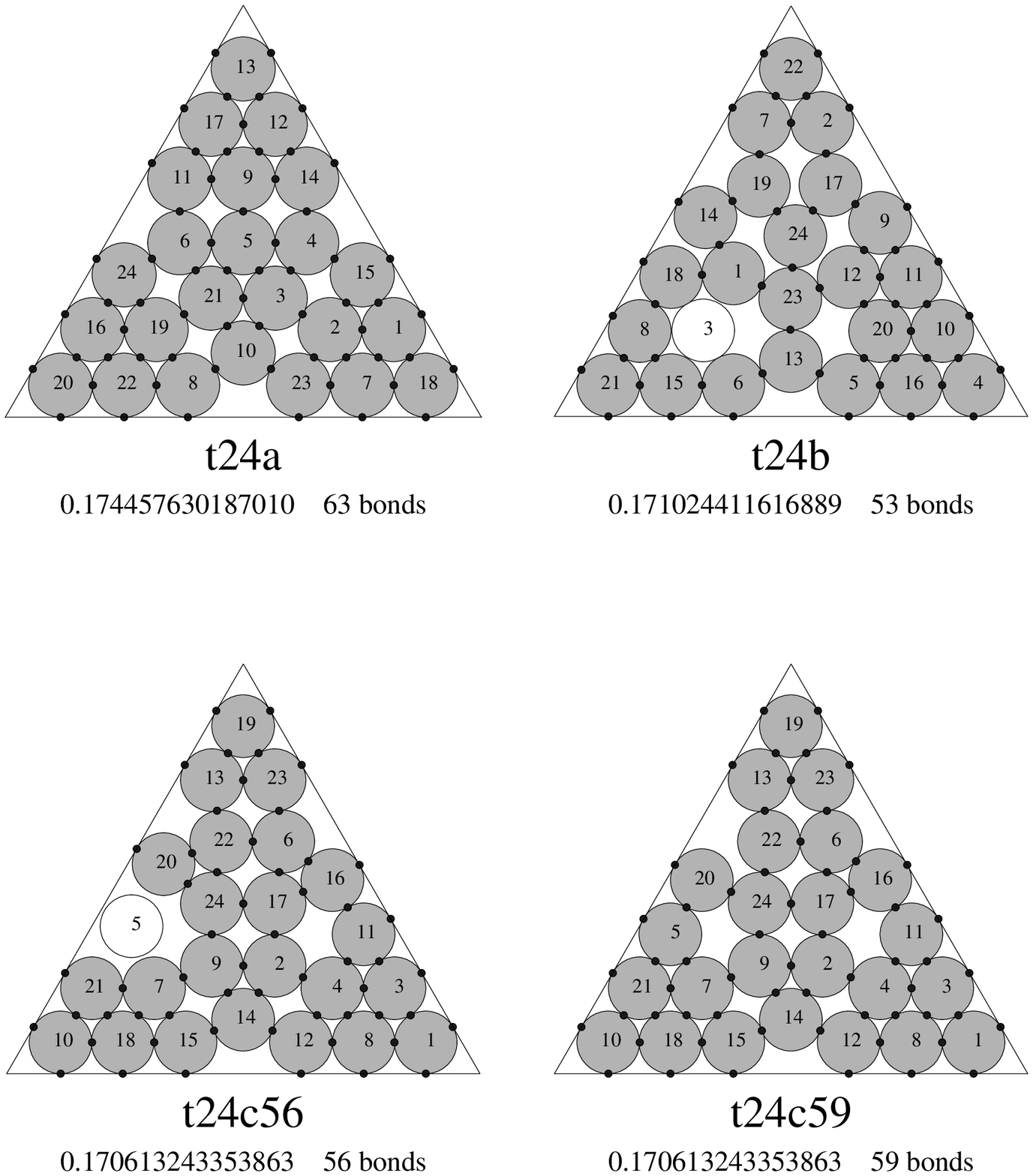,width=7in}}

\caption{The best (t24a), the next-best (t24b), and the third-best (t24c56, t24c59) packings of 24 disks.}
\label{t24}
\end{figure}
\begin{figure}[htb]
\centerline{\psfig{file=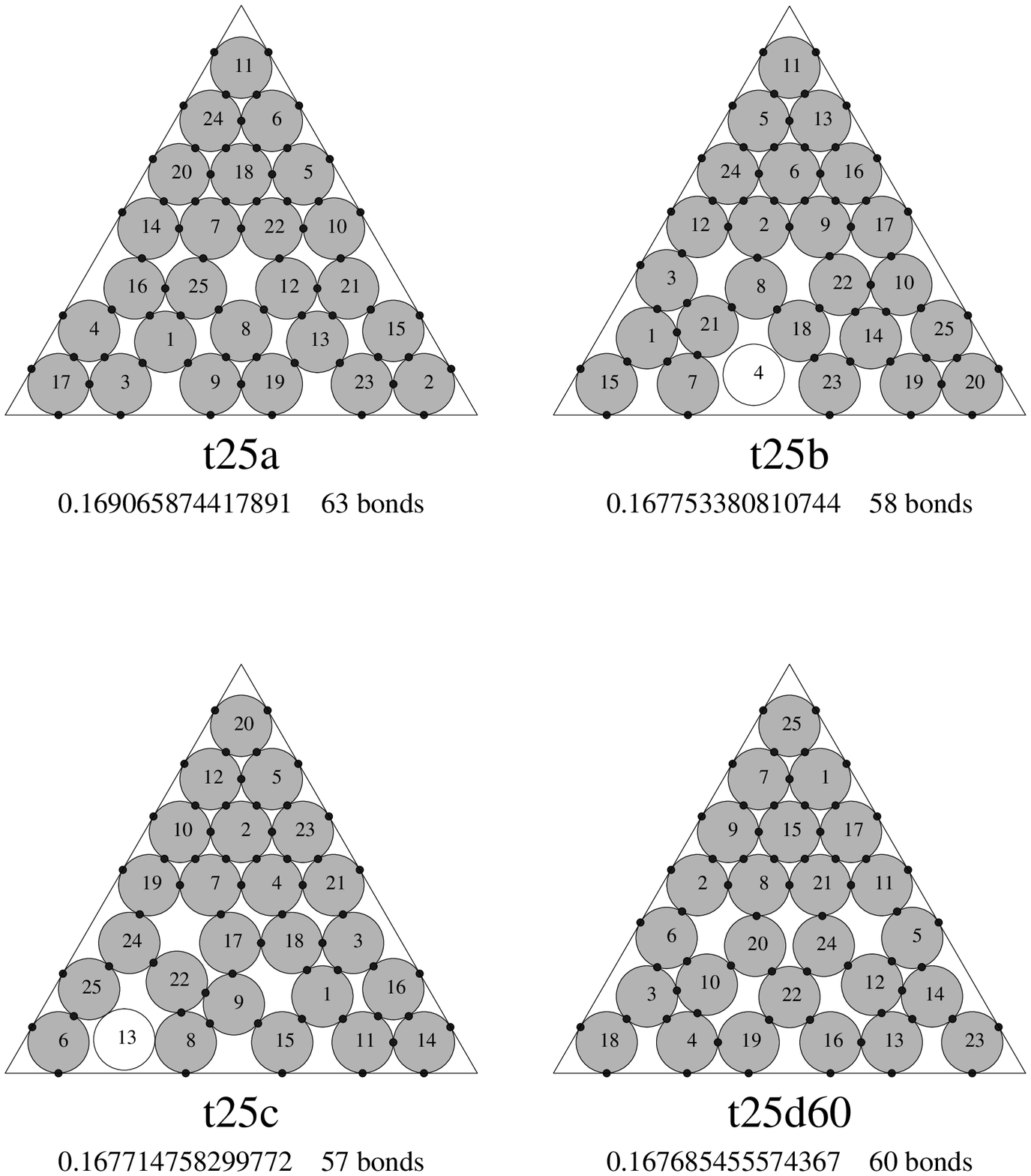,width=7in}}

\caption{The best (t25a), the next-best (t25b), the third-best (t25c), and a fourth-best (t25d60) packings of 25 disks.}
\label{t25}
\end{figure}
\begin{figure}[htb]
\centerline{\psfig{file=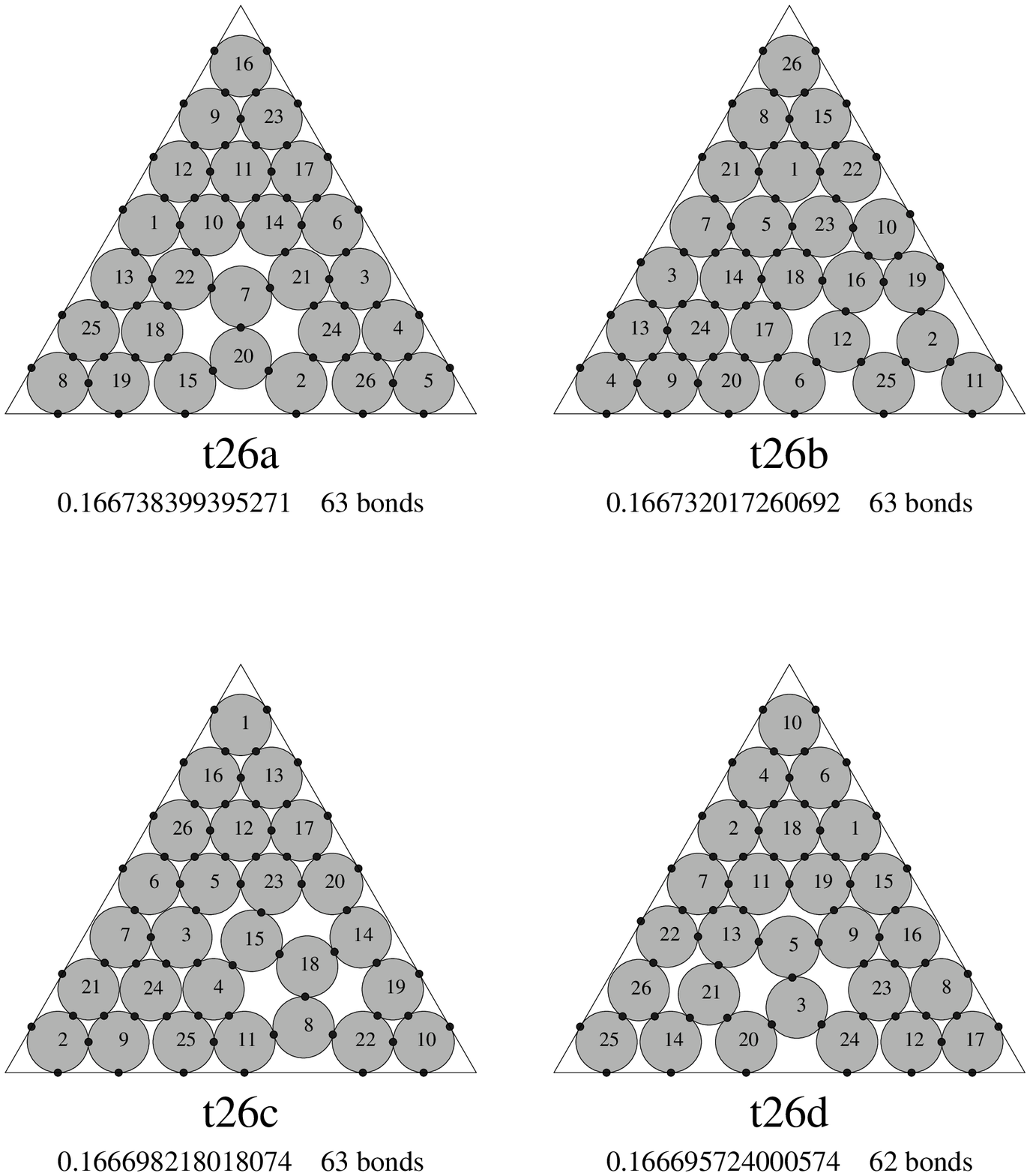,width=7in}}

\caption{The best (t26a), the next-best (t26b), the third-best (t26c), and the fourth-best (t26d) packings of 26 disks.}
\label{t26}
\end{figure}
\begin{figure}[htb]
\centerline{\psfig{file=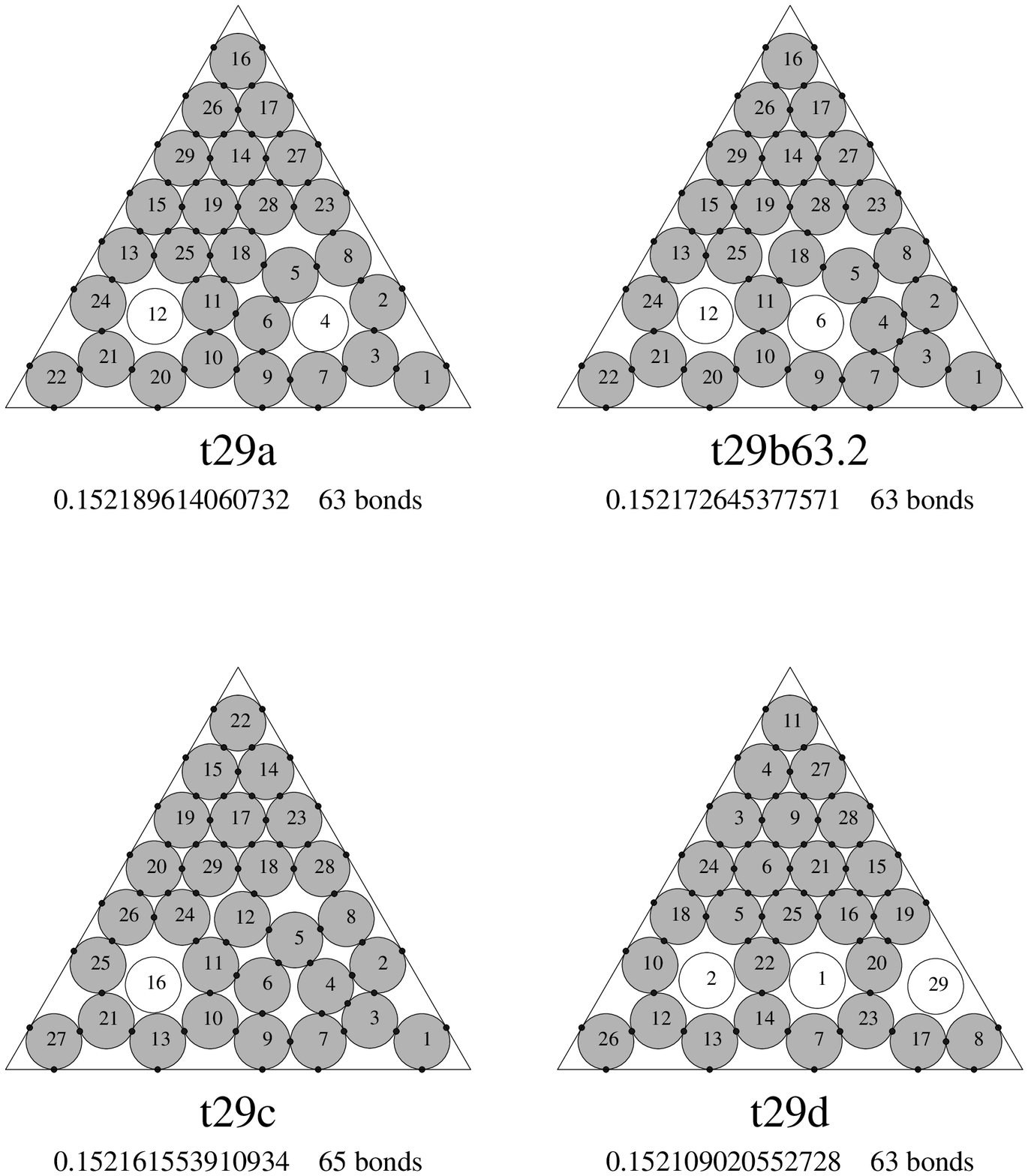,width=7in}}

\caption{The best (t29a), a next-best (t29b63.2), the third-best (t29c), and the fourth-best (t29d) packings of 29 disks.}
\label{t29}
\end{figure}
\begin{figure}[htb]
\centerline{\psfig{file=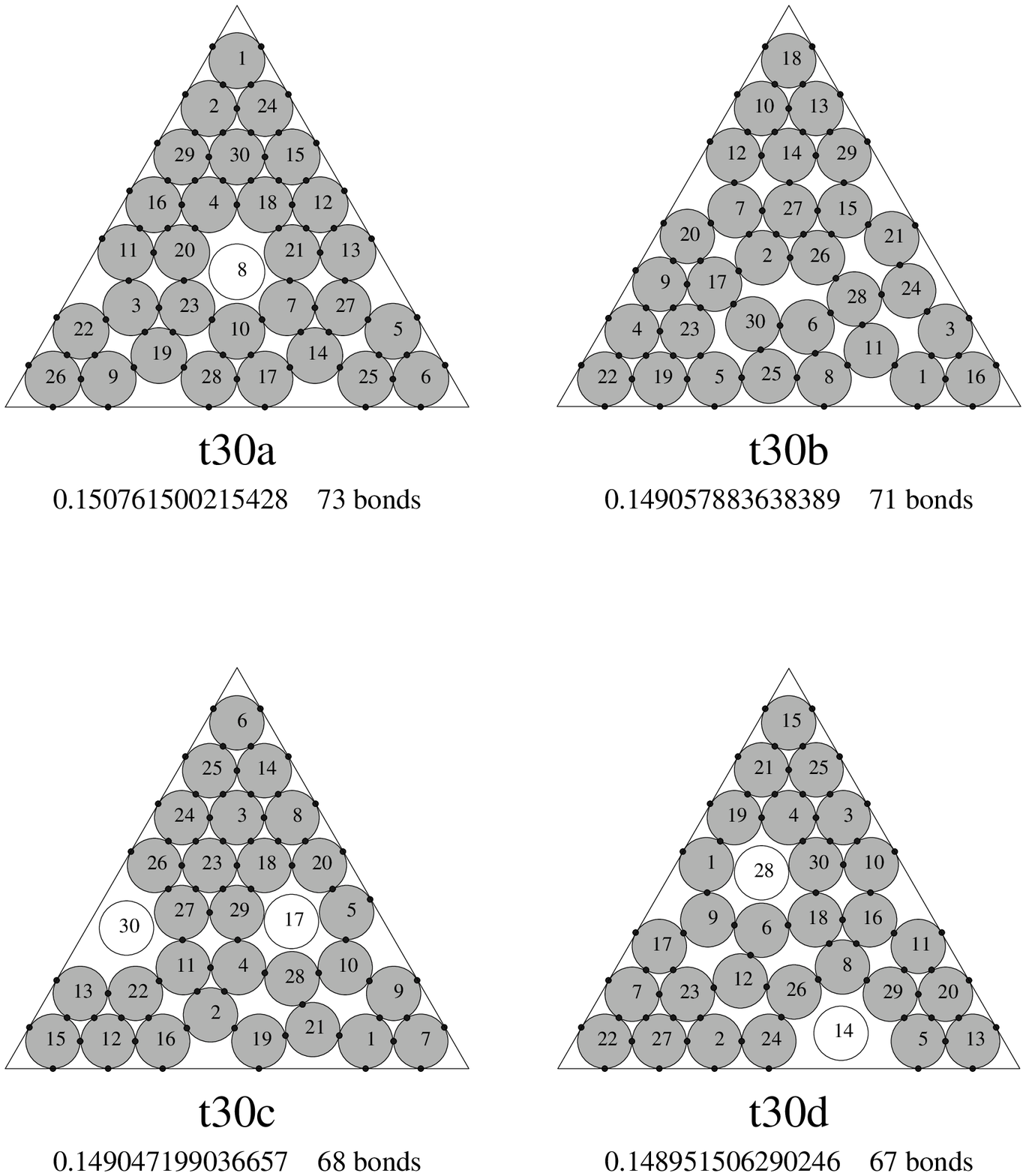,width=7in}}

\caption{The best (t30a), the next-best (t30b), the third-best (t30c), and the fourth-best (t30d) packings of 30 disks.}
\label{t30}
\end{figure}
\begin{figure}[htb]
\centerline{\psfig{file=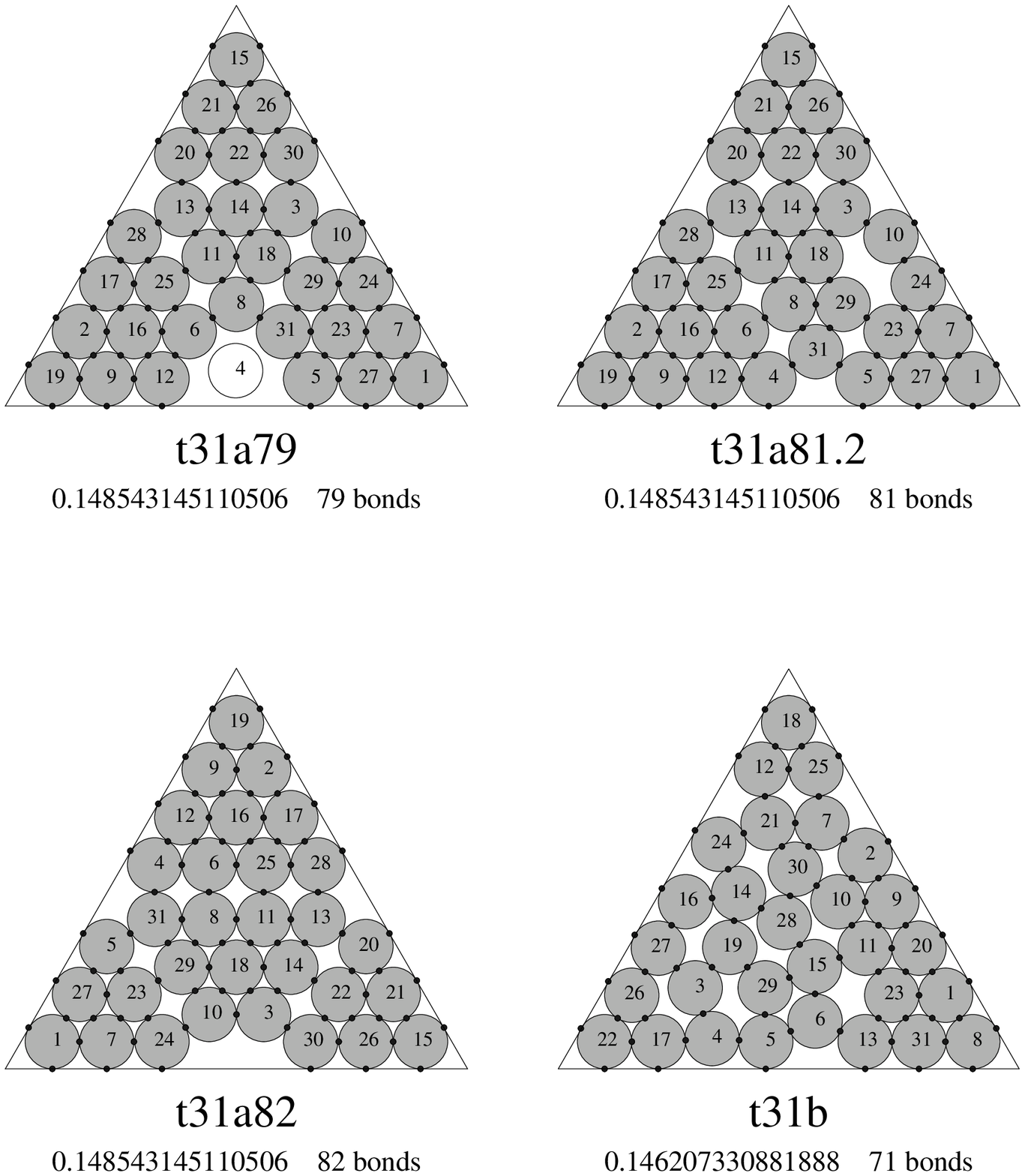,width=7in}}

\caption{The best (t31a79, t31a81.2, t31a82) and the next-best (t31b) packings of 31 disks.}
\label{t31}
\end{figure}
\begin{figure}[htb]
\centerline{\psfig{file=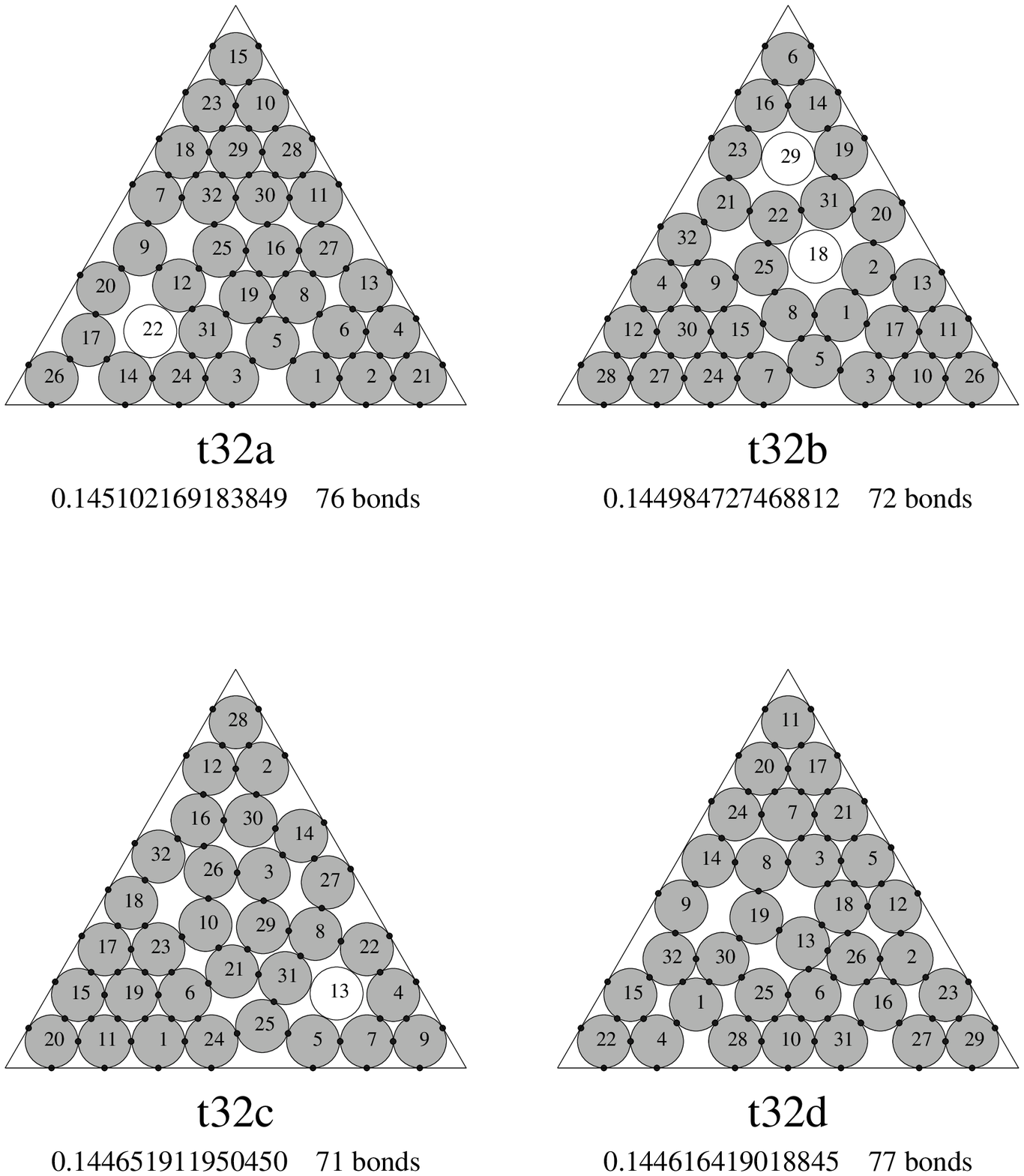,width=7in}}

\caption{The best (t32a), the next-best (t32b), the third-best (t32c), and the fourth-best (t32d) packings of 32 disks.}
\label{t32}
\end{figure}
\begin{figure}[htb]
\centerline{\psfig{file=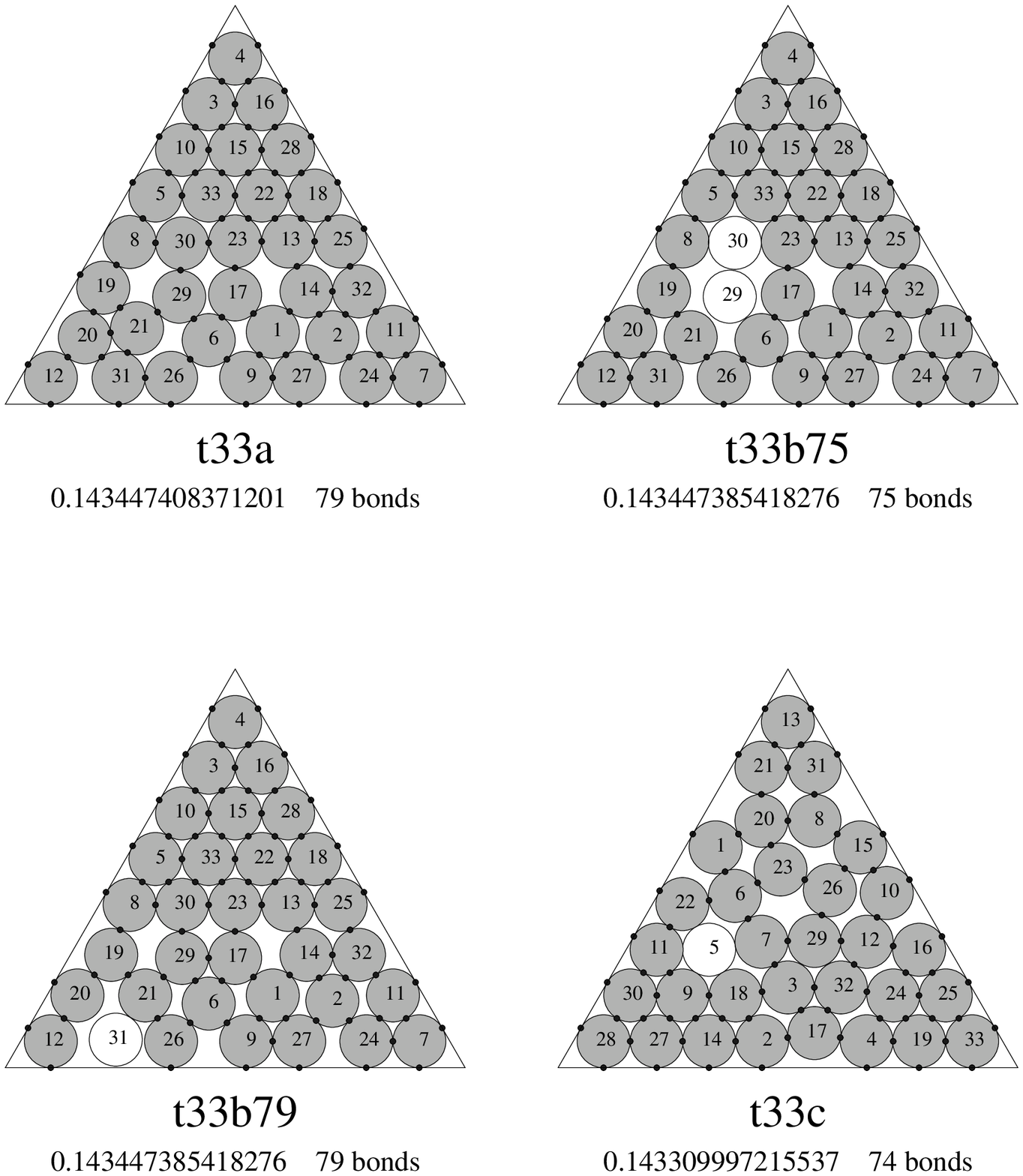,width=7in}}

\caption{The best (t33a), the next-best (t33b75, t33b79), and the third-best (t33c) packings of 33 disks.}
\label{t33}
\end{figure}
\begin{figure}[htb]
\centerline{\psfig{file=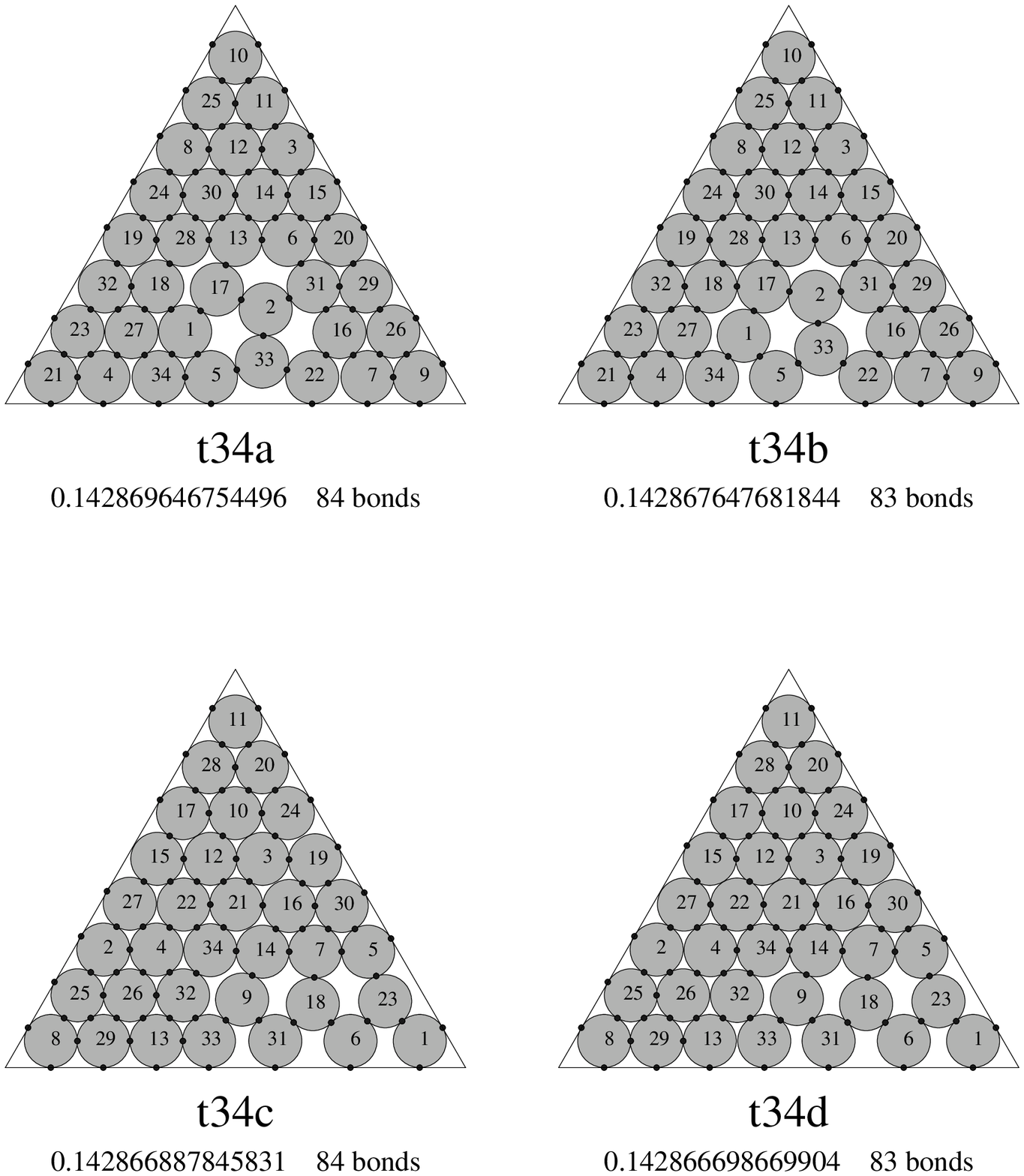,width=7in}}

\caption{The best (t34a), the next-best (t34b), the third-best (t34c), and the fourth-best (t34d82) packings of 34 disks.}
\label{t34}
\end{figure}

\clearpage
\section{Conjectures for individual packings}
\hspace*{\parindent}
Each packing diagram we give can imply several different statements:
\begin{itemize}
\item[(I)]
There exists a valid configuration of nonoverlapping disks with all pairwise distances marked by bonds equal to zero, and those not marked by bonds strictly positive, and with disk diameter equal to the indicated value with a relative error of less than $10^{-14}$.
\item[(II)]
The configuration is rigid:
no disk or set of disks except for rattlers can be continuously displaced from the indicated positions without overlaps.
\item[(III)]
The configurations are correctly ranked.
That is, the a-packing really is optimal, the b-packing is second best, etc.
\end{itemize}

We believe (I) and (II) are correct.
As to (III), we hope the statement is correct with respect to the a-packings.
In other words, we believe these are the optimal packings.
We are less confident for the lower ranked packings.
For example, if someone finds a new packing in between our
c- and d-packings, we will not be astounded.
We provide these mainly for comparison purposes, and to serve as benchmarks for other packing algorithms.

It would also not be surprising 
to discover a nonisomorphic packing to one we have
presented which has exactly the same disk diameter and the same number of bonds
(e.g., as in t17b42ns and t17b42s in Fig.~\ref{t17}).
\begin{figure}[htb]
\centerline{\psfig{file=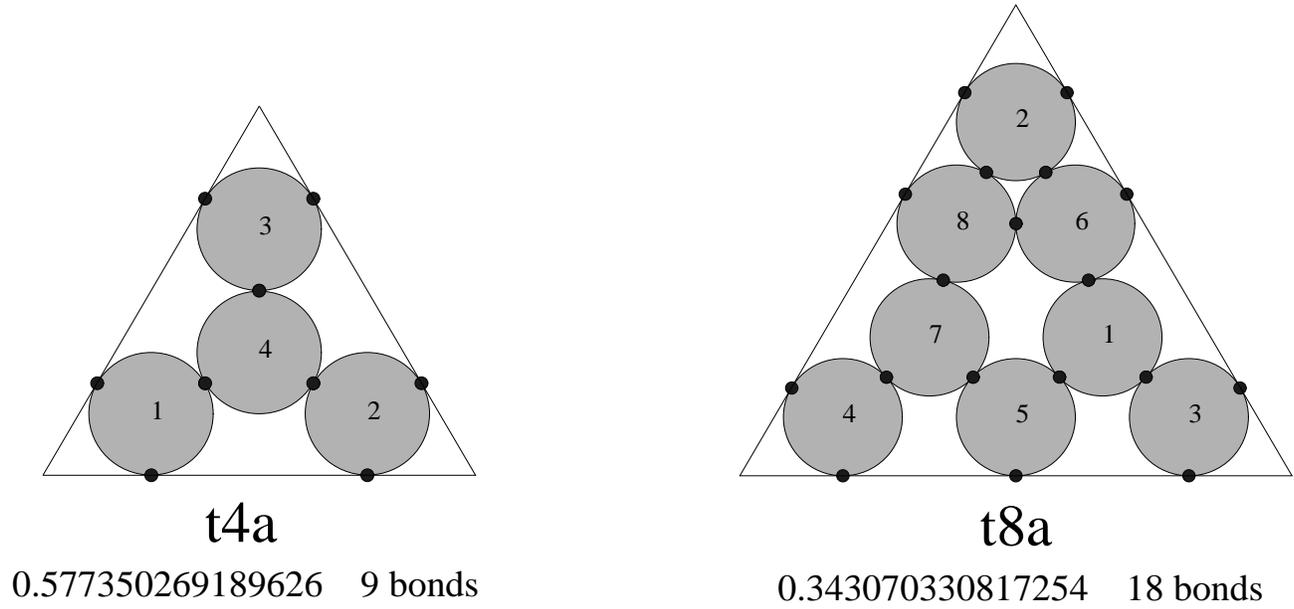,width=7in}}
\caption{The best packing of 4 disks (t4a) and 8 disks (t8a).}
\label{t4-8}
\end{figure}
\section{Conjectures for infinite classes}
\begin{figure}[htb]
\centerline{\psfig{file=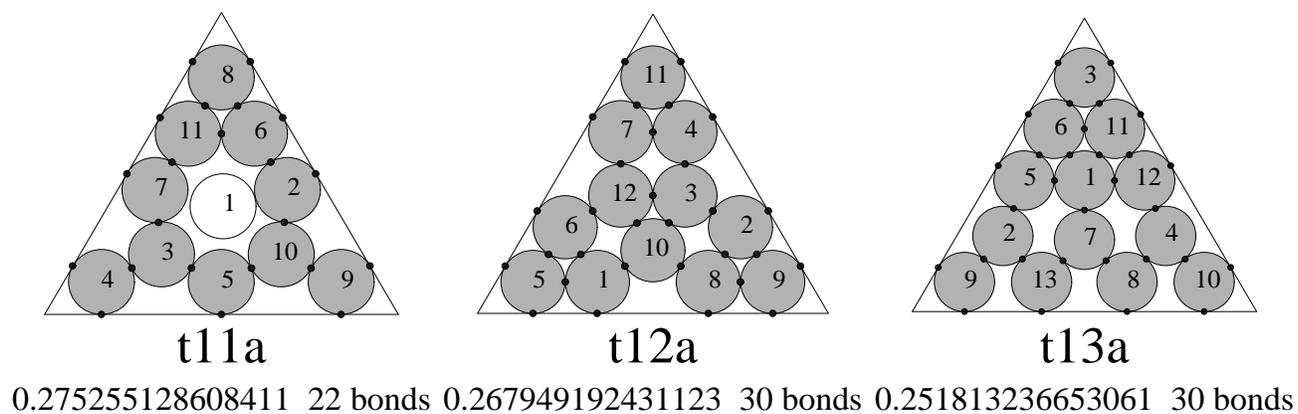,width=7in}}

\caption{The best packings of 11 disks (t11a), of 12 disks (t12a), and of 13 disks (t13a).}
\label{t11-13}
\end{figure}

\begin{figure}[htb]
\centerline{\psfig{file=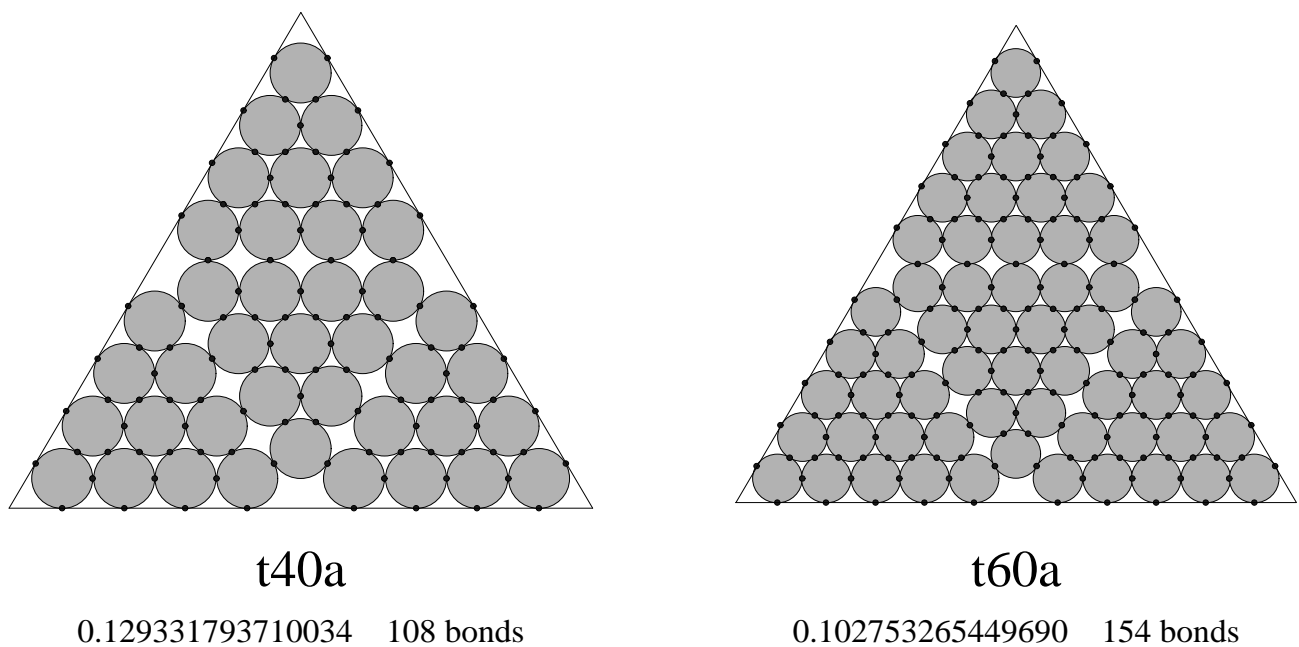,width=7in}}

\caption{The best packings of 40 disks (t40a) and 60 disks (t60a).}
\label{t40-60}
\end{figure}

\hspace*{\parindent}
Dense packings in an equilateral triangle seem to ``prefer'' to form blocks of dense triangles and arrangements that are nearly so.
In this section we describe seven infinite classes 
where we think we have found the optimal packings.
Each class has its individual pattern of the optimal packings,
which is different from the patterns for other classes.
However, since they are the result of the particular packings we found,
which themselves are only conjectured to be optimal,
then the general conjectures have even less reliability.
We still think they might serve as useful organizers for the maze of published dense packings.
\paragraph{4$\Delta$(k).}
The best packing (we found) of $24 = 4 \Delta (3)$ disks 
in Fig.~\ref{t24} consists of four triangles, each with $\Delta (3) =6$ disks.
The best packings of 12 disks (in \cite{M1}) in Fig.~\ref{t11-13}, and
even 4 disks in Fig.~\ref{t4-8} have the same form.
The packings we obtained while experimenting 
with $40 = 4 \Delta (4)$ and $60 = 4 \Delta (5)$ disks 
(see Fig.~\ref{t40-60}), and also with 
$84 = 4 \Delta (6)$ and $112 = 4 \Delta (7)$ disks
have the same structure as well.

A simple analysis of the patterns obtained implies that
$d(4 \Delta (k)) = \frac{1}{ 2k - 2 + \sqrt{3}}$.

If we fit members of class 4$\Delta$(.) within the
boundaries of the triangular periods,
i.e., among members of the class $\Delta$(.), then
every other period 
has exactly one $n$ of the form 4$\Delta(k)$ 
lying almost exactly at the middle of the period.

\begin{figure}[htb]
\centerline{\psfig{file=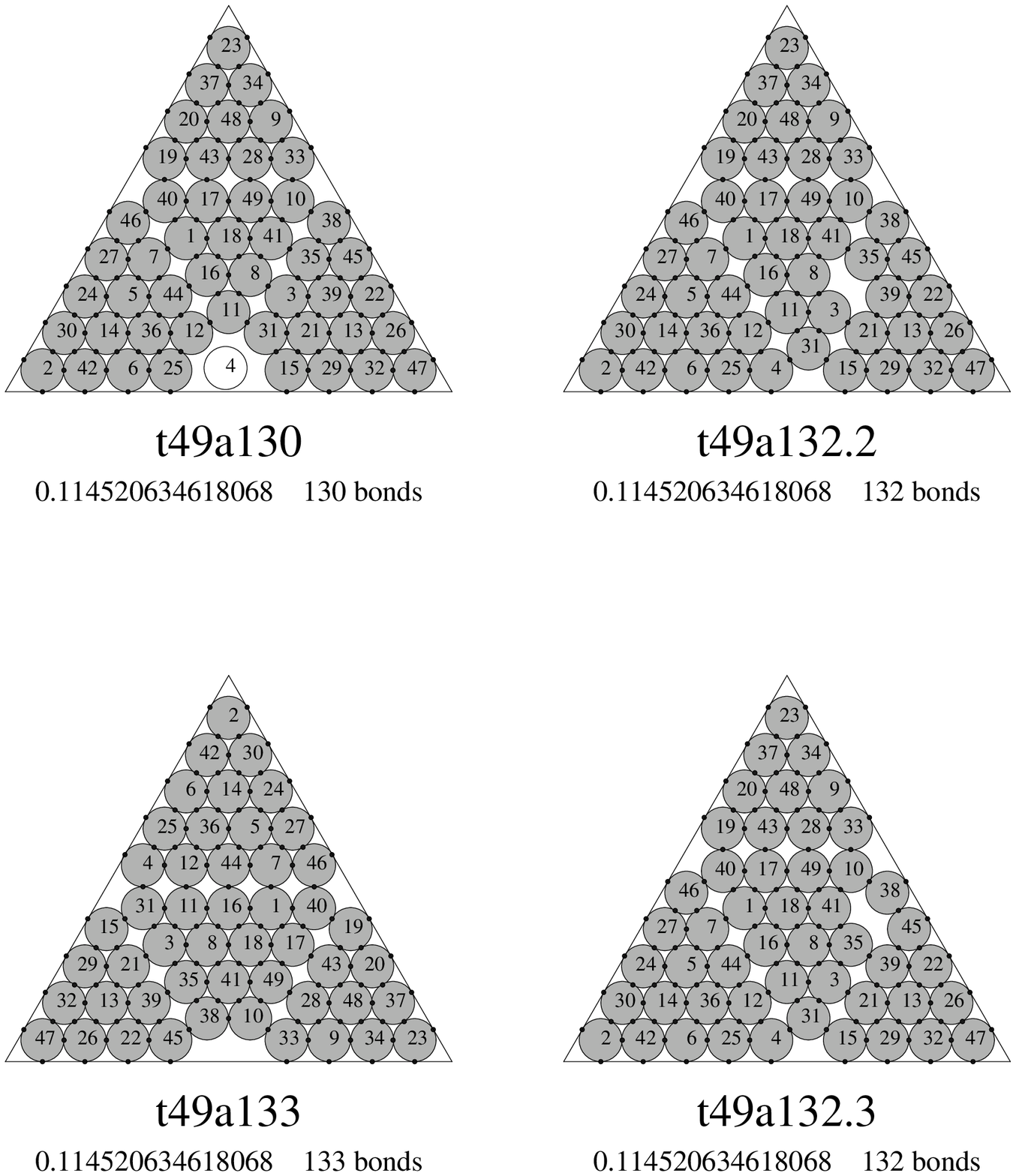,width=7in}}

\caption{Four (out of the five existing) best packings of 49 disks.}
\label{t49}
\end{figure}

\paragraph{${\bf 2} \mbox{\protect\boldmath $\Delta$} {\bf (k+1) + 2} \mbox{\protect\boldmath $\Delta$} {\bf (k) -1}$.}
For each $k$ there are $k+1$ distinct best packings:
two for 7 disks (Fig.~\ref{t7}), three for 17 disks (Fig.~\ref{t17}),
four for 31 disks (three of these four are shown in Fig.~\ref{t31}),
and five for 49 disks (four of these five are shown in Fig.~\ref{t49}).
\clearpage
\begin{figure}[htb]
\centerline{\psfig{file=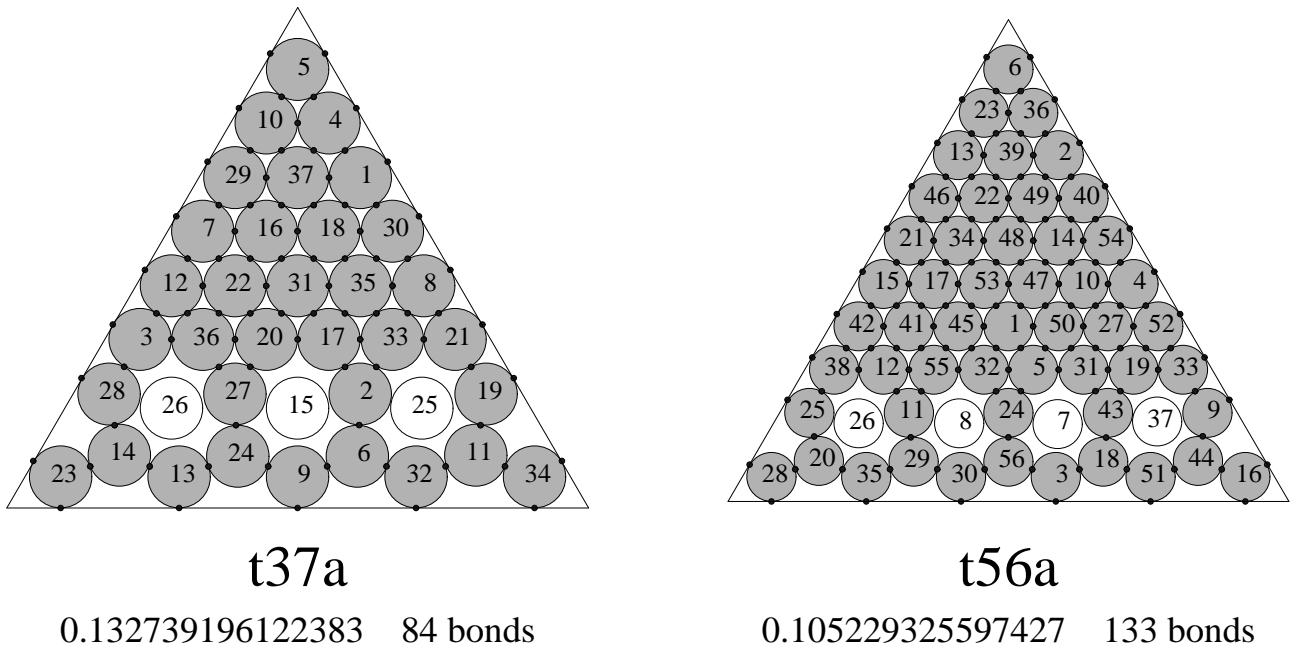,width=7in}}

\caption{The best packings of 37 disks (t37a) and 56 disks (t56a).}
\label{t37-56}
\end{figure}
\begin{figure}[htb]
\centerline{\psfig{file=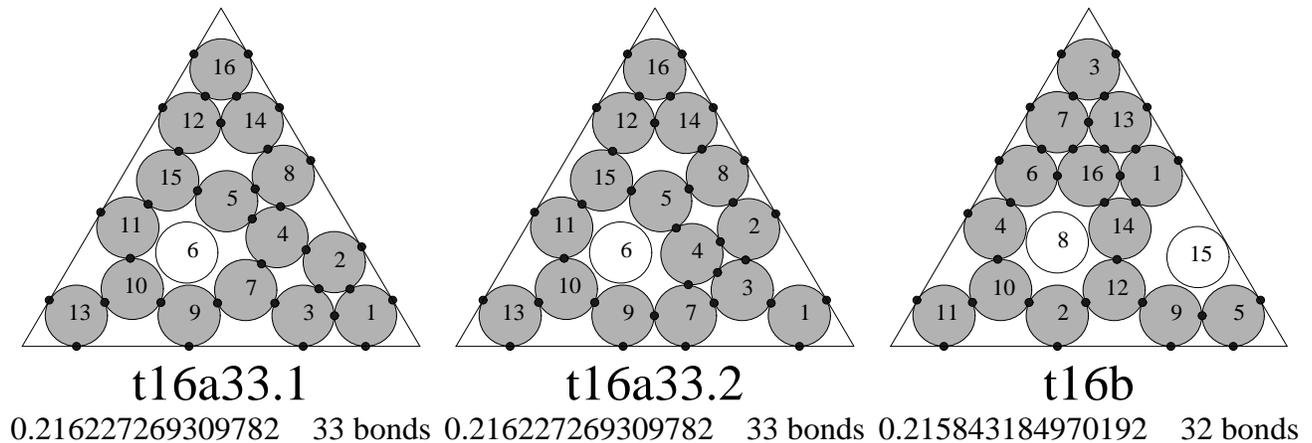,width=7in}}

\caption{The best (t16a33.1, t16a33.2) and the next-best (t16b) packings of 16 disks.}
\label{t16}
\end{figure}
\begin{figure}[htb]
\centerline{\psfig{file=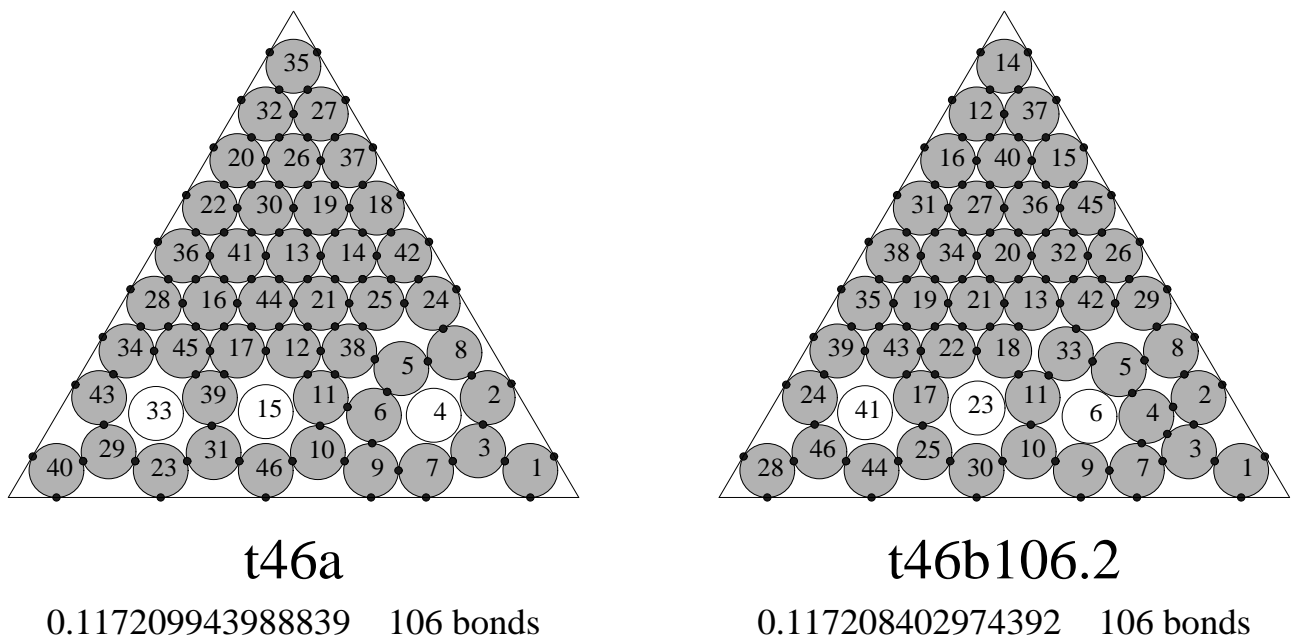,width=7in}}

\caption{The best (t46a) and a next-best (t46b106.2) packings of 46 disks.}
\label{t46}
\end{figure}
\begin{figure}[htb]
\centerline{\psfig{file=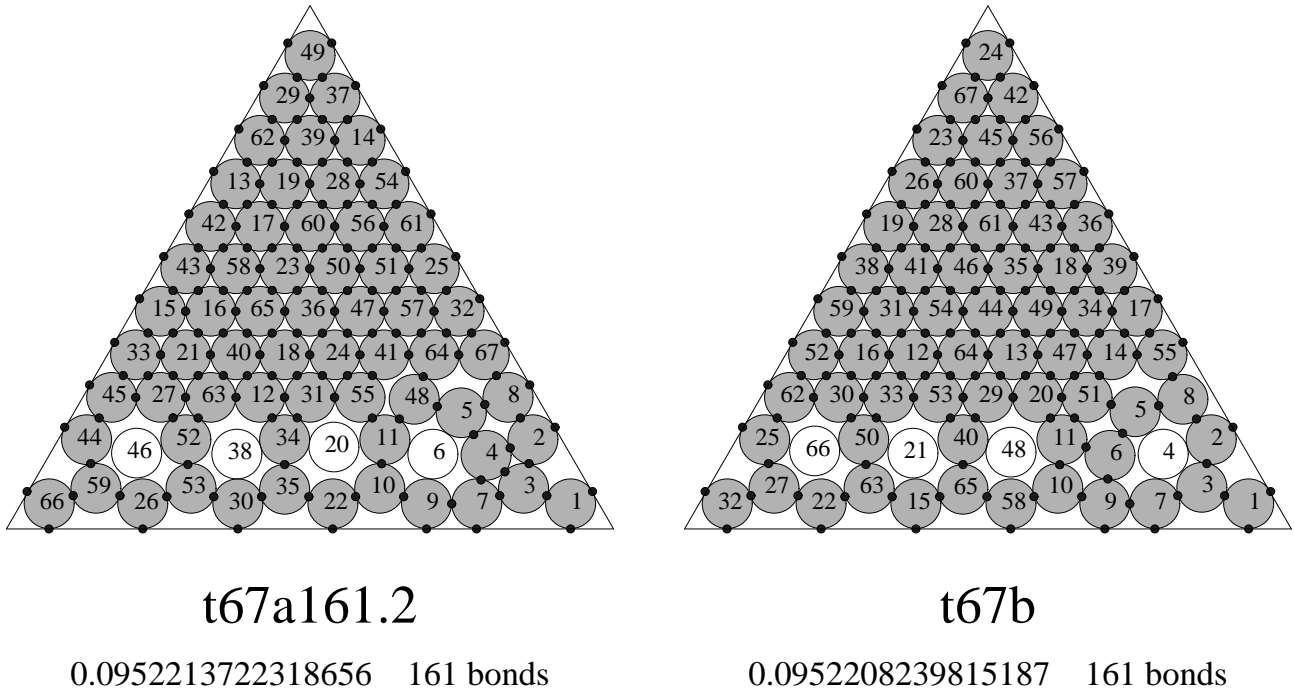,width=7in}}

\caption{A best (t67a161.2) and the next-best (t67b) packings of 67 disks.}
\label{t67}
\end{figure}
\begin{figure}[htb]
\centerline{\psfig{file=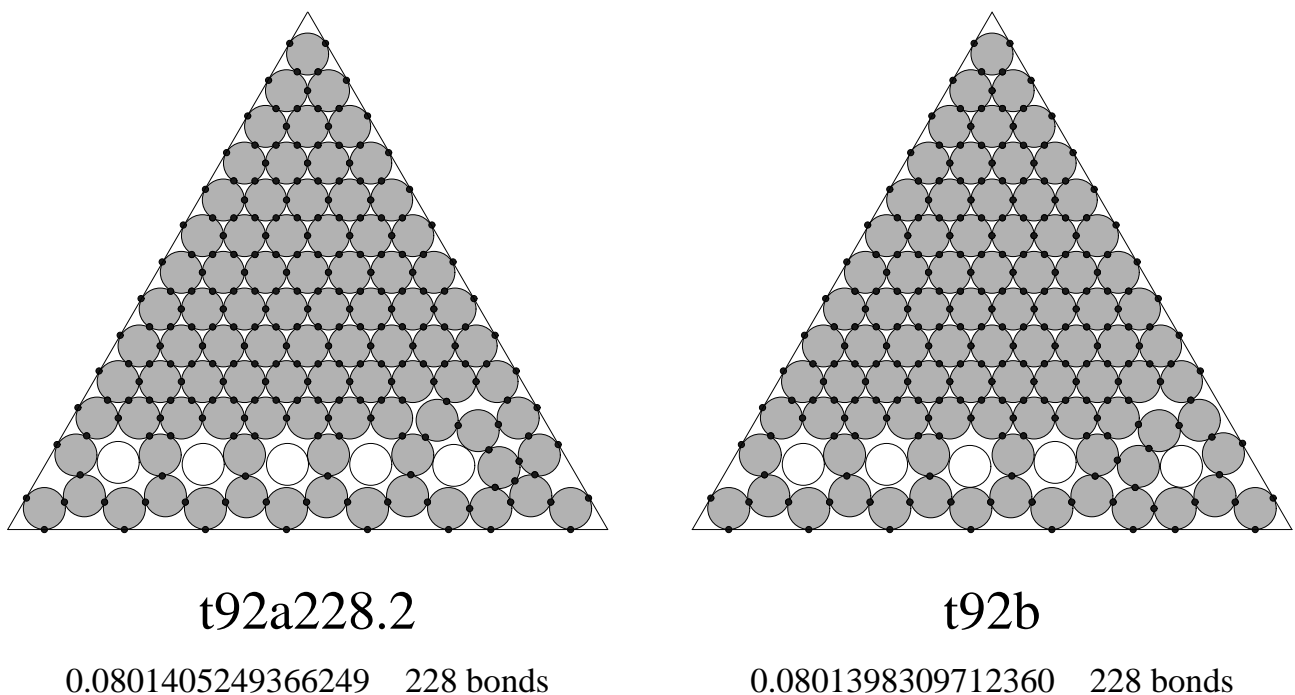,width=7in}}

\caption{A best (t92a228.2) and the next-best (t92b) packings of 92 disks.}
\label{t92}
\end{figure}
\begin{figure}[htb]
\centerline{\psfig{file=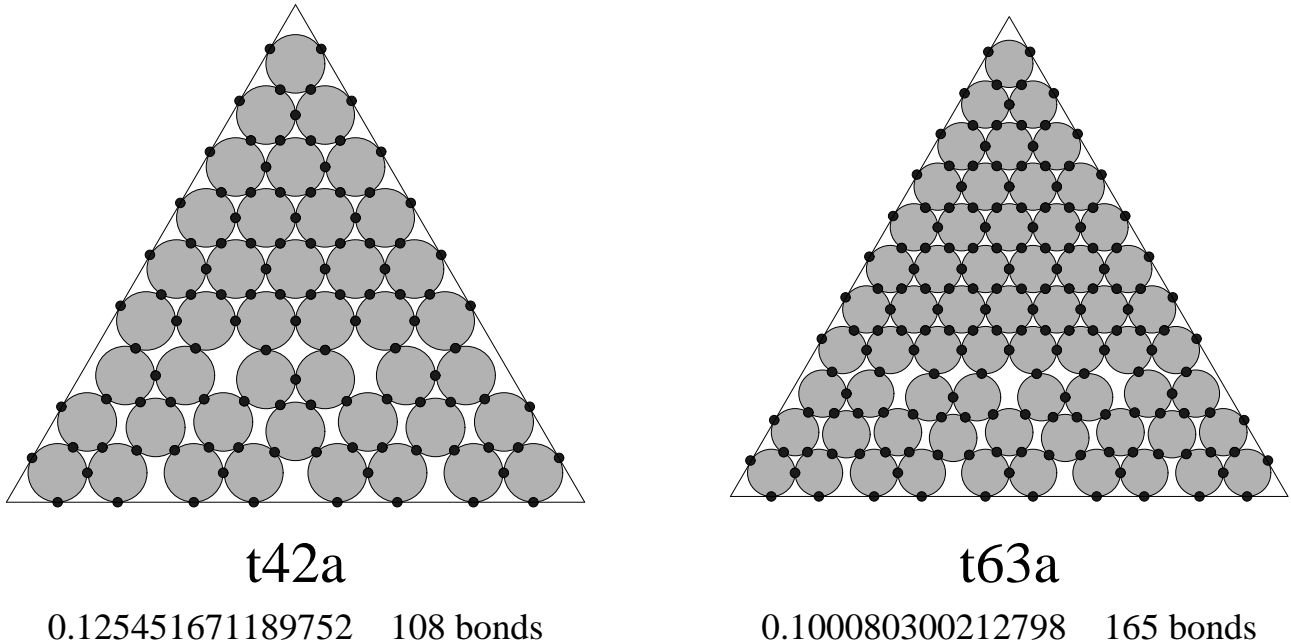,width=7in}}

\caption{The best packings of 42 disks (t42a) and of 63 disks (t63a).}
\label{t42-63}
\end{figure}
\begin{figure}[htb]
\centerline{\psfig{file=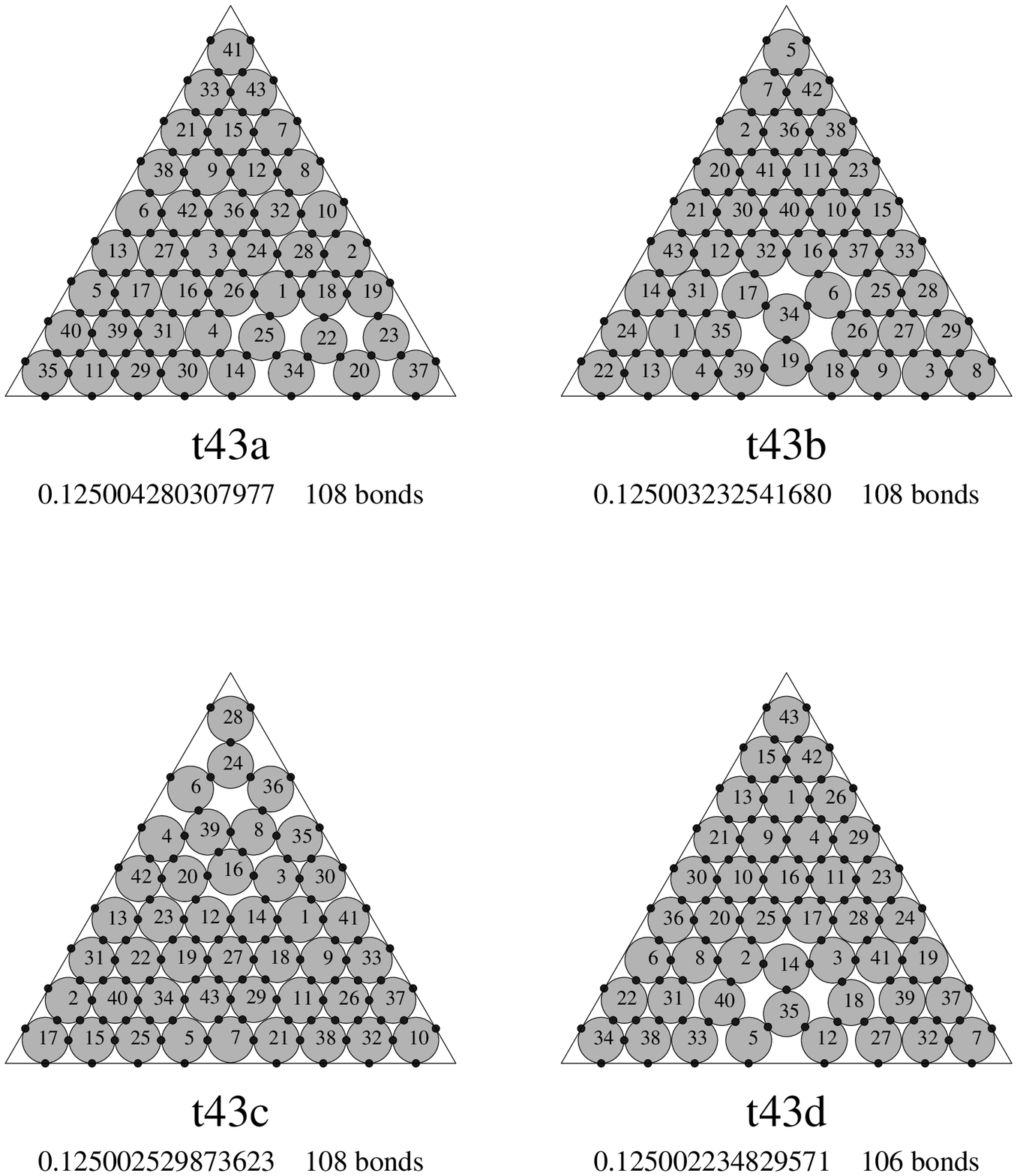,width=7in}}

\caption{The best (t43a), the next-best (t43b), the third-best (t43c) and the fourth-best (t43d) packings of 43 disks.}
\label{t43}
\end{figure}
\begin{figure}[htb]
\centerline{\psfig{file=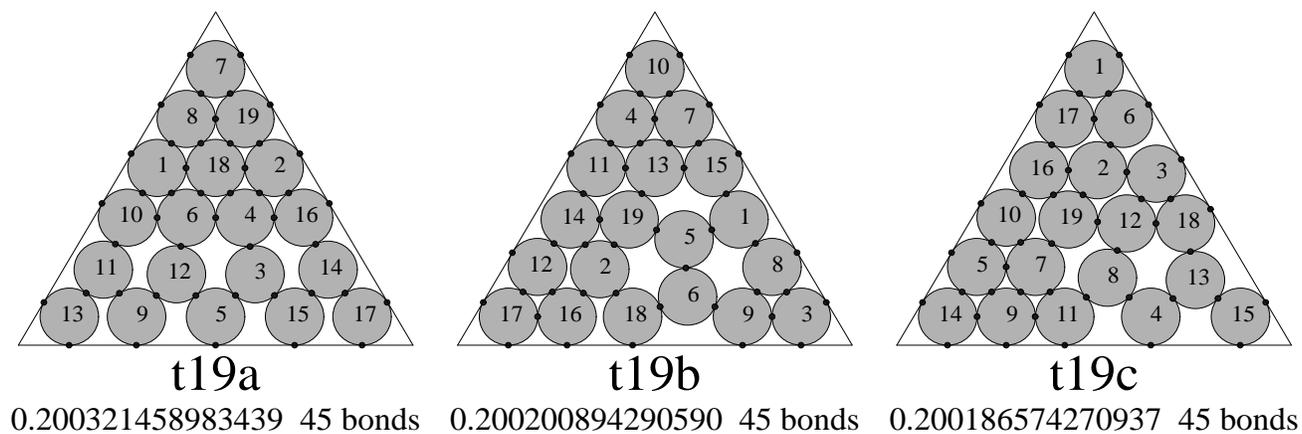,width=7in}}

\caption{The best (t19a), the next best (t19b), and the third best (t19c) packings of 19 disks.}
\label{t19}
\end{figure}
\begin{figure}[htb]
\centerline{\psfig{file=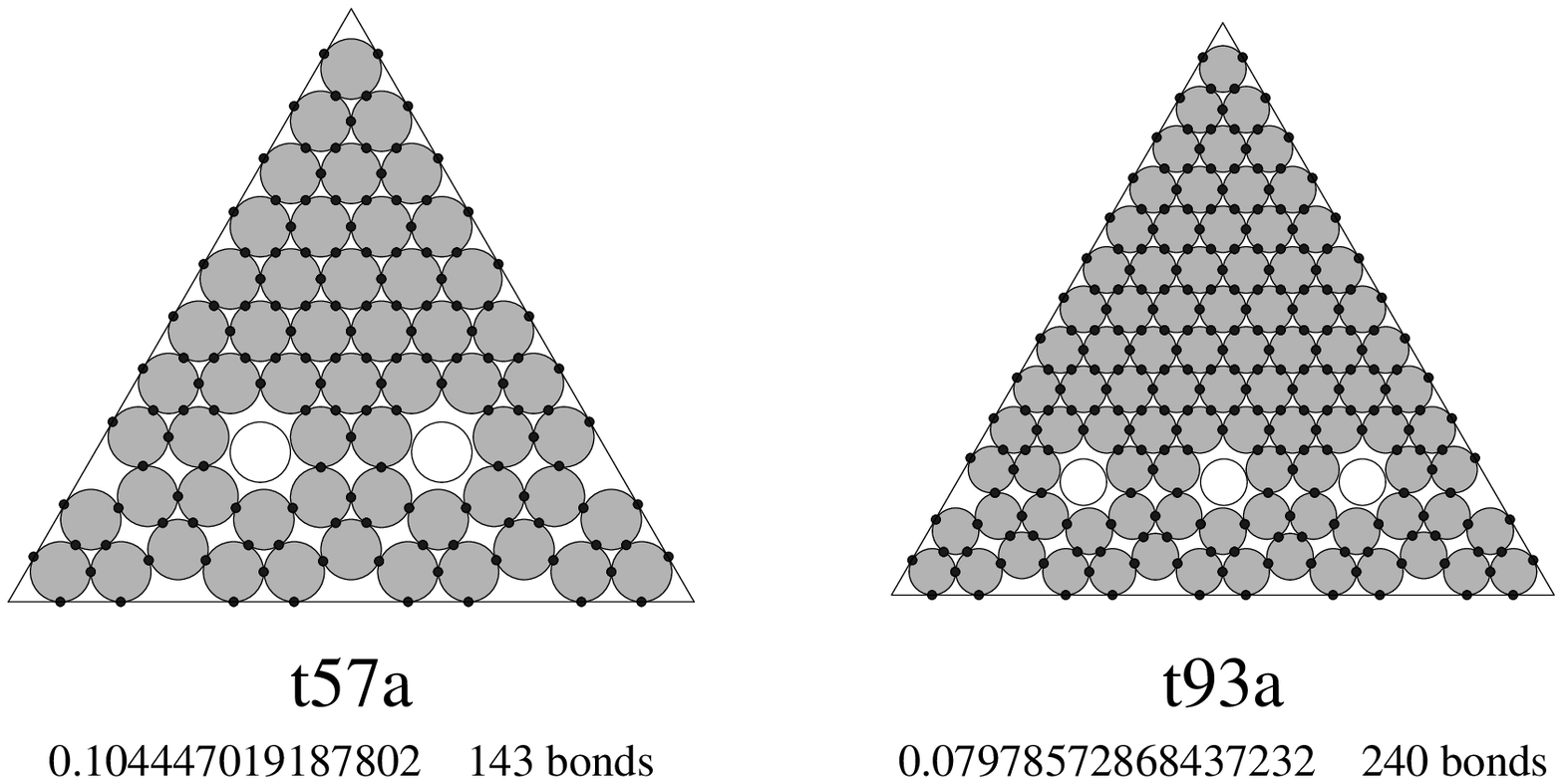,width=7in}}

\caption{The best packings of 57 disks (t57a) and 93 disks (t93a).}
\label{t57-93}
\end{figure}

There are three different numbers of bonds in these packings;
the smallest number of bonds is in the packing with a rattler,
then $k-2$ packings each of which has a ``cavity'' 
and the same number of bonds,
and finally one more packing without a rattler or a cavity with the
largest number of bonds.
The four triangles, two small and two large,
that illustrate the expression $2 \Delta (k+1) + 2 \Delta (k)-1$,
can be seen in the packings with a rattler
(t17a40,  t31a79, t49a130).
The two larger triangles are defective: both coalesce
a corner disk (disk 6 in t17a40, disk 4 in t31a79, disk 4 in  t49a130).
Packing t31a81.2 is obtained from t31a79 by the
left larger triangle acquiring its
corner disk 4
and pushing disks 31 and 29 from the left side of the other large triangle 
down into the cavity formed .
If we push only disk 31, we obtain a packing t31a81.1 (not shown).
If we push all three disks 31, 29, 10 into the cavity
(and rotate the resulting structure to
recover the symmetry with respect to the vertical axis),
we obtain the fourth best packing t31a82.

A simple analysis of the patterns obtained implies that
$d(2 \Delta (k+1) + 2 \Delta (k) -1) = \frac{1}{ 2k - 1 + \sqrt{3}}$.
Our experiments with $71 = 2 \Delta(6)+2 \Delta(5) -1$ disks
produced the same patterns for the best packing. 

If we fit the members of class $2\Delta (k+1)+2\Delta (k)-1$ into the
boundaries of the triangular periods, as we did for the class $4 \Delta(k)$,
we find that every other period
has exactly one $n$ of the form $2\Delta(k+1)+2\Delta (k)-1$
which lies almost exactly in the center of the period.
Thus, classes $2\Delta (k+1)+2\Delta (k)-1$ and
$4 \Delta(k)$ are ``parity-complementary'' to each other. 
Beginning with the second triangular period 3 to 6, each 
period has exactly one term of one or the other class;
even periods contain terms of the class $4 \Delta(k)$
and odd periods contain terms of the class $2\Delta(k+1)+2\Delta (k)-1$.

The pattern of one of the parity-complementary class pair can be obtained
from the pattern of the other class by a simple transformation.
For example if we eliminate the 
eight bottom disks 2, 42, 6, 25, 15, 29, 32, 47, and the rattler 4
in t49a130 (Fig.~\ref{t49}),
we obtain the pattern of t40a (Fig.~\ref{t40-60}).

\paragraph{$\mbox{\protect\boldmath $\Delta$} {\bf (2k)+1}$.}
As we noted earlier,
$\Delta (m)$
densely packed disks in an equilateral triangle 
form a perfect hexagonal lattice with 
$m$ disks on a side.
When $m=2k$ is even,
the structure with $\Delta (2k)+1$ disks adjusts itself to
one extra disk as follows.
The top $2k-1$ rows remain packed hexagonally, and the bottom 
row ripples  to accommodate $2k+1$ disks instead of $2k$.
In this ripple of the bottom row,
the 1st, 3rd, ... $(2k+1)$th disk beginning from the left corner
remain attached to the bottom,
while the 2nd, 4th, ... $2k$th disk rise
and  attach themselves
to 2nd, 4th, ... $2k$th  disks respectively,
of the row above.
Thus, $k-1$ rigid cages are formed. 
The $k-1$ disks of the row above
to which no disk is attached from below
fall off into these cages and become rattlers.
The first seven terms of the class $\Delta (2k)+1$ are:
t4a ($k = 1$, no rattlers since $k-1 = 0$, see Fig.~\ref{t4-8}),
t11a (constructed in \cite{M1} with one rattler; see Fig.~\ref{t11-13}), 
t22a (see Fig.~\ref{t22}) with two rattlers, t37 and t56 (Fig~\ref{t37-56})
with 3 and 4 rattlers, respectively,
t79a (194 bonds, 5 rattlers, $d(79)=0.0871159038791759$),
and t106a (267 bonds, 6 rattlers, $d(106)=0.0742982999063026$).
We do not reproduce the diagrams here for the latter two packings;
their patterns are identical to the class description given above.
\paragraph{$\mbox{\protect\boldmath $\Delta$} {\bf (2k+1)+1}$.}
When $m=2k+1$, $k=1,2...$, the odd parity of $m$ 
causes a more complex adjustment to the extra disk.
The bottom row ripples in a non-symmetric way;
the ripple creates $k$ cages for rattlers
and a cavity; see t29a (Fig.~\ref{t29}) for $k = 3$.

Notice that packing t29b63.2 (Fig.~\ref{t29}) has almost
the same structure as t29a, 
except for the cage that consists of disks 2, 3, 7, 9, 6, 5, and 8,
is depressed and disk 4 in t29b63.2 is not a rattler,
and a nonrattler 6 in t29a becomes a rattler in t29b63.2.
The same two modifications exist for $k = 4$ (i.e., $n = 46$),
and the modification with the depressed cage, t46b106.2, is again inferior
(Fig.~\ref{t46}).
Beginning with $k = 5$ (i.e., $n = 67$),
while both modifications exist, they exchange their roles:
the depressed one becomes the best, t67a161.2, 
while the other one becomes
the inferior one, t67b (Fig.~\ref{t67}).
For example, t92a228.2 is the modification with the depressed
cage, while t92b is the other one (Fig.~\ref{t92}).
The same pattern is displayed by packing t121a307.2 
(307 bonds, 6 rattlers, $d(121)=0.0691630188894699$),
for which we omit the diagram here.

Labels t29b63.2, t46b106.2, t67a161.2, t92a228.2,  and t121a307.2
have the suffix 2 in them
because there exist equivalent packings 
t29b63.1, t46b106.1, t67a161.1, t92a228.1, and t121a307.1, respectively.
The latter differ from the former in the placement of only 4 disks.
An easy way to explain this is to look at the second term of the class
$n=\Delta(2k+1)+1$ 
for $n = 16$ disks (Fig.~\ref{t16}).
In this case both modifications
exist and both deliver the optimum, t16.a33.1 and t16.a33.2.
They differ in the placement of disks 2, 3, 4, and 7.

The side rattler disk 15 in t16b can be considered a precursor
for the side rattler 29 in t29d.
The same side-rattler pattern was observed in lower ranked
packings for
$k = 4$ ($n = 46$), $k = 5$ ($n = 67$), $k = 6$ ($n = 92$), 
and $k = 7$ ($n = 121$).

Classes $\Delta(2k)+1$ and $\Delta(2k+1)+1$ are a parity complementary
pair, similar to the pair of classes 
$4\Delta(k)$ and $2\Delta(k+1)+2\Delta(k)-1$
considered above.
\paragraph{$\mbox{\protect\boldmath $\Delta$} {\bf (k+2)-2}$.}
While the optimal packings of $\Delta (k+2) -1$ disks are always
(apparently) perfectly hexagonal with a single disk removed, 
the removal of two disks from a hexagonal arrangement is never optimal.
The first two terms 
t4a and t8a (Fig.~\ref{t4-8}) suggest no common pattern.
Looking at the next case t13a (Fig.~\ref{t11-13}) 
suggests the pattern for $k \ge 3$ 
of a packed triangle $\Delta (k)$ at the top,
supported by two sparse rows of disks, each of which lacks a disk
compared to what would be there in a perfect hexagonal packing.
The top-$\Delta (k)$-plus-two-sparse-rows packing
indeed exists for any $k \ge 3$ and is rigid.
In particular, the pattern appears again in the best packings 
for $k=4$ (t19a in Fig.~\ref{t19}).

However, the optimality of this pattern
does not continue for $k > 4$,
as can be seen
in Fig.~\ref{t26} where the best packing t26a 
has a different pattern.
The top-$\Delta (k)$-plus-two-sparse-rows pattern is not even among 
the top four packings for $n=26$.
The pattern for t26a and t26b persists for the next term 
(see t34a and t34b in Fig.~\ref{t34})
but then roles become reversed for 43 disks 
(see t43a and t43b in Fig.~\ref{t43}).
Will the pattern of packing t43a remain optimal for larger values of $k$?
Unfortunately, our algorithm fails to obtain 
stable packings for 53 (or larger values of $\Delta (k+2) -2$) disks.
\paragraph{$\mbox{\protect\boldmath $\Delta$} {\bf (3k+1) +2 =} \mbox{\protect\boldmath $\Delta$} {\bf (3k-1) + (2k+1)} \mbox{\protect\boldmath $\Delta$} {\bf (2)}$.}
$30=\Delta (5)+5\Delta (2)$ and packing t30a (Fig.~\ref{t30}) 
can be viewed as a $\Delta (5)$ triangle on top 
from which disk 8 fell off and became a rattler, supported by five
triangles $\Delta (2)$ from below;
t30a is the second term of the class.
To produce this structure for the next value of $k$
we add two triangles $\Delta (2)$ on the bottom and 
three more layers of disks to
enlarge the top triangle to become a $\Delta (8)$ (and, in general,
this procedure is repeated for larger values of $k$).
Indeed our experiments with $n=57$ ($k=3$) and 
$n=93$ ($k=4$) did produce this structure in the best
packings (see t57a and t93a in Fig.~\ref{t57-93}).
For $k=1$ we have $n=12$, a degenerate case with no rattlers.
The number of rattlers in this packing for general $k$ is $k-1$.
\paragraph{$\mbox{\protect\boldmath $\Delta$} {\bf (2k+3) -3 =} \mbox{\protect\boldmath $\Delta$} {\bf (2k) + (2k+1)} \mbox{\protect\boldmath $\Delta$} {\bf (2)}$.}
The first term is packing t12a (Fig.~\ref{t11-13}) 
which also belongs to class $4 \Delta(.)$.
The second term is packing t25a (Fig.~\ref{t25})
which can be viewed as a $\Delta (2k)$ triangle on top 
supported by $2k+1$ alternating $\Delta (2)$ triangles below.
This pattern of the second term, whose description also fits the first term,
is more apparent in the third and fourth terms
(see t42a and t63a shown in Fig.~\ref{t42-63}).
According to our experiments the next terms t88a and t117a
do not continue this pattern
but give way to patterns which are somewhat
similar to those of the class $\Delta(k+2)-2$ considered above.
(We omit diagrams for both t88a and t117a.)
\section{How good are the packings?}
\hspace*{\parindent}
Let us compare the values of the best packings
with the only bound currently available, namely one
based on the inequality of Oler \cite{O}.
This inequality has the following form 
(see \cite{FG} for a simple proof):
Let $K$ be a compact, 
convex subset of $\EE^2$ with area $A(K)$ and perimeter $P(K)$.
If $p(K)$ denotes the maximum number of points that can be placed
in $K$ so that any pair has mutual distance at least 1,
then the following inequality holds:
\beql{eq1}
p(K) \le \frac{2}{\sqrt{3}} A(K) + \frac{1}{2}
P(K) +1 ~.
\eeq
Inverting (\ref{eq1}) and applying it with $p(K) =n$ and $K$ being an
equilateral triangle of side length $L(n)$, we obtain
\beql{eq2}
L(n) \ge \frac{1}{2} ( -3 + \sqrt{8n+1}) : = t (n) ~.
\eeq

\begin{figure}[H]
\centerline{\psfig{file=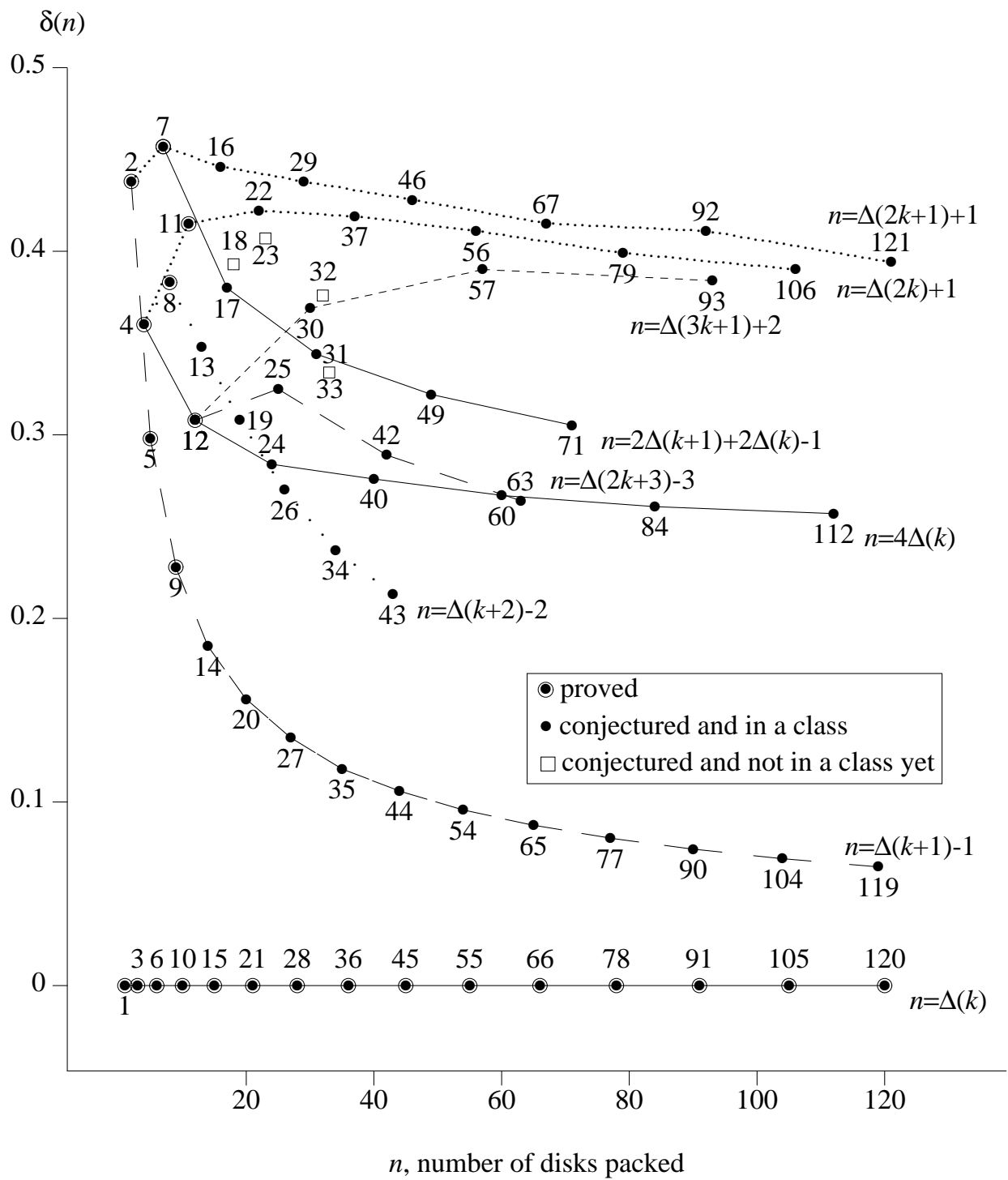,width=6.7in}}

\caption{Discrepancy between side length of the triangle and its lower bound for different $n$.}
\label{fig61}
\end{figure}

In Fig~\ref{fig61} we plot the difference
$$\delta (n) : = L(n) - t(n)$$
versus $n$ for selected values of $n \le 121$.
A dot with a circle around it indicates that the corresponding 
value has been proved to be optimal;
a dot without a surrounding circle or an open square
indicates that the value is only conjectured to be optimal.
For up to 37 disks there were only four values which
did not fall into one of our infinite classes, namely, $n$ = 18, 23, 32, 
and 33; those are indicated by open squares.
The values which have been associated
with classes are connected by lines,
with a distinct type of line for each class.

We should point out that for each $n$, the value of the largest disk
diameter $d(n)$ and the value of $L(n)$ are reciprocally related, i.e.,
$d(n) L(n) =1$.
Thus, $L( \frac{k(k+1)}{2} ) =k-1$ for $k \ge 1$.
If our conjecture for
$n= \Delta (k) +1 = \frac{k(k+1)}{2} +1$ is correct then it would follow that
$$\lim_{k \to \infty} \delta \left( \frac{k(k+1)}{2} +1 \right) \le
\frac{2}{\sqrt{3}} -1 = 0.1547 ~\ldots ~.
$$

Using the explicit values for $d(4\Delta(k))$ and
$d(2\Delta(k+1)+2\Delta(k)-1)$ given in Section~5, 
we have $L(4\Delta(k)) = 2k -2 +\sqrt{3}$  and 
$L(2\Delta(k+1)+2\Delta(k)-1) =2k -1 +\sqrt{3}$, from which it
follows that both $\lim_{k \to \infty} \delta (4\Delta (k))$ and
$\lim_{k \to \infty} \delta (2\Delta(k+1)+2\Delta(k)-1)$
are at most $\sqrt{3} - 3/2 = 0.2321 ~\ldots ~$
with the distance between $k$th term and the limit
being of the order of $1/k$.

We conjecture that for each of the classes
$n=\Delta(2k)+1$, $\Delta(2k+1)+1$, and $\Delta(3k+1)+2$ 
the value of $\delta (n)$ is bounded away from zero. 
In fact, we believe that for any fixed $c > 0$,
the value of $\delta (\Delta(k)+c)$ is bounded away from zero. 
In other words, packings of $\Delta(k)$ disks are so tight
that any attempt to accommodate even one additional disk
noticeably worsens the packing in that $L(n)$ increases by
at least some positive amount independent of $n$. 
In this sense the packings for the class $\Delta(k)$ are ``tight''.

On the other hand, we believe that after any fixed positive number 
of disks are added
to $\Delta(k)$ disks, any other fixed number of disks can be added
without substantial ``damage'' to $\delta(n)$ (asymptotically).
Thus, for example, it would seem 
that if $\lim_{k \to \infty} \delta(n(k))$ exists for each
of the class
$n(k)=\Delta(2k)+1$, $\Delta(2k+1)+1$, and $\Delta(3k+1)+2$
(and the limits probably do exist), then all three limits are equal.
In this sense the classes $\Delta(2k)+1$, $\Delta(2k+1)+1$, 
and $\Delta(3k+1)+2$ are ``loose''.

Similarly, the classes $n = \Delta(k+2)-2$ and $\Delta(2k+3)-3$ are 
``loose'' in the sense that we can add one disk to the best
(conjectured) packing
without noticeable change of $\delta(n)$ for sufficiently large $n$.
This follows from the fact that
$\lim_{k \to \infty} \delta(\Delta(k)-p) \to 0$
for $p$ fixed.
The latter limit is obvious by noticing
that a lower-bounding packing for $\Delta(k)-p$ disks
when $k$ is sufficiently large
is simply the densest packing of $\Delta(k)$ disks
with $p$ disks removed;
for such a packing, $\delta$ is asymptotically 0.

Formally,
we say that an infinite class of packings of $n$ disks,
$n=n(1), n(2),...n(k),...$,
is {\em loose}, if
$\lim_{k \to \infty}$[$\delta(n(k)+1) - \delta(n(k))$] $=0$.
Because we believe this limit exists for any class
we consider, each class has to be either tight or loose.

We further conjecture that the classes $4\Delta(k)$
and $2\Delta(k+1)+2\Delta(k)-1$ are tight, similar to the class $\Delta(k)$.

Are there other tight classes?
Here is our argument in favor
of the existence of a countable
{\em infinity} of distinct tight classes.
We believe that if each densest packing of
$n$ disks for $n = n(1), n(2),$...$n(k)$,.., 
consists of a {\em fixed} number, say $r$,
of densely packed triangles $\Delta (.)$,
then $\delta(n(k)+1) - \delta(n(k))$ is bounded away from zero
as $k$ goes to infinity.

Thus, we have the following task: for each $r$ from some infinite
set find a sequence $n(1), n(2),$...$n(k)$,..,
so that the densest packing of $n(k)$ disks for all $k$
has the ``same pattern'' and consists of exactly $r$ densely
packed triangles.
Note that {\em a priori} we are not able to define
what the ``same pattern'' is (and hence we have
no formal definition of what a ``class'' is);
but after producing a class the pattern is usually clear.

\begin{figure}[H]
\centerline{\psfig{file=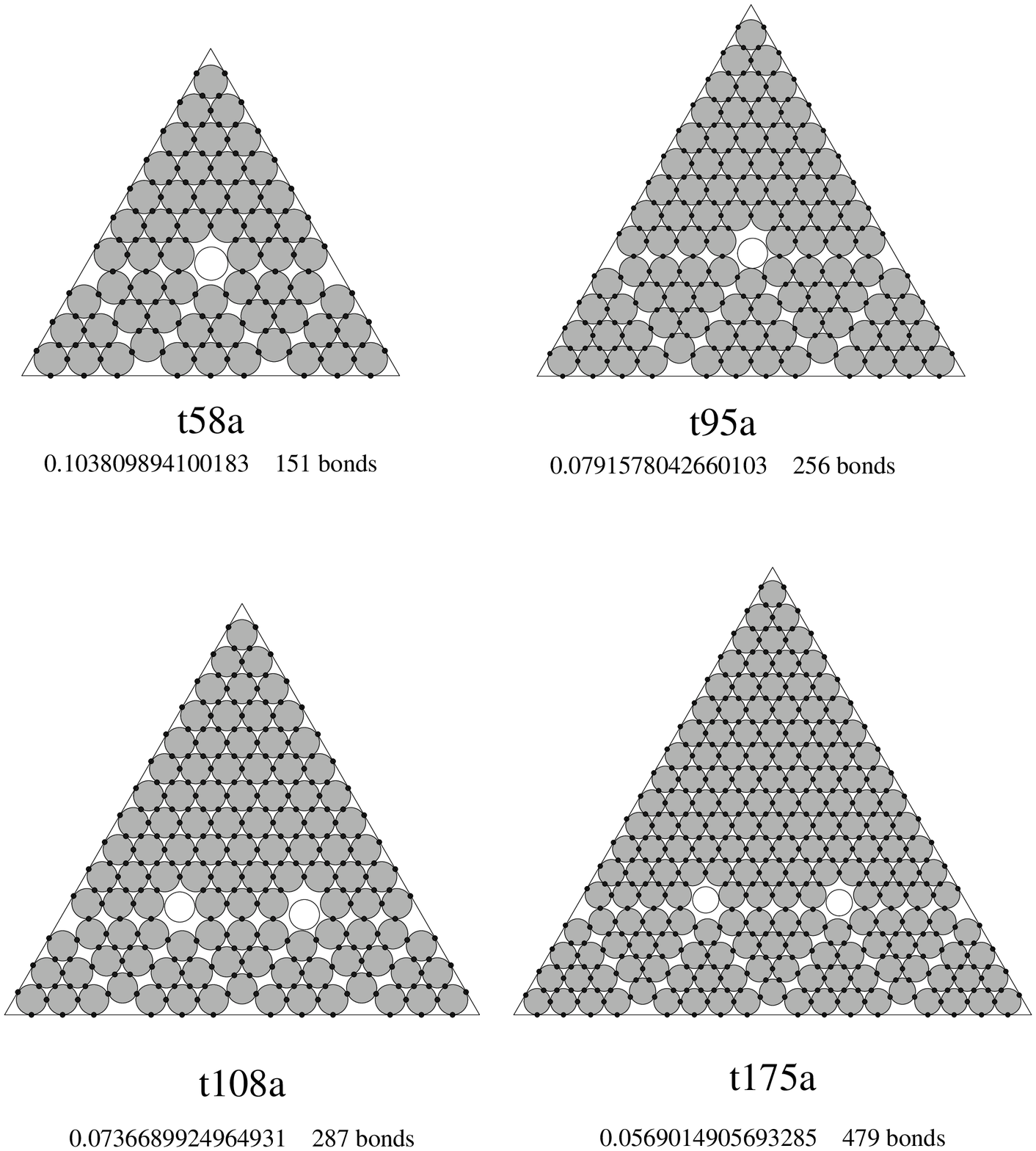,width=7.1in}}

\caption{The best packings of 58 disks (t58a), 95 disks (t95a), 108 disks (t108a), and 175 disks (t175a).}
\label{t58-175}
\end{figure}

Let us consider
a two-parameter family of numbers $n= n_p (k),~ p=1,2...,~ k=1,2...$,
of the form
\beql{eq3}
n_p (k) = \Delta((k+1)(p+1)-2) + k~=~\Delta((k+1)p-1) + (2p+1) \Delta(k)~.
\eeq
Two equal expressions for $n_p (k)$ are given in (3).
The second expression suggests $r = 2(p+1)$ triangles.
If we take, for example, $k=p=2$, we get $n_2(2)=30$ disks
and $r=6$.
The conjectured 
t30a (Fig~\ref{t30}) indeed consists of 6 densely packed
triangles, if we attach rattler 8 to the top triangle
$\Delta(5)$.
Take now 
$58=n_2 (3)$. 
The pattern of our experimental packing t58a
(Fig~\ref{t58-175}) looks like t30a (Fig~\ref{t30}) with 6 triangles again
(with the rattler attached to the top triangle).

Thus, just as t30a is a member of the class $\Delta(3k-1) + (2k+1) \Delta(2)$
for $k=2$,
t58a is (perhaps) a member of the class $\Delta(4k-1) + (2k+1) \Delta(3)$
for $k=2$.  
The pattern of the $k$th packing of this class 
is composed of $2(k+1)$ densely packed triangles.
Is this class tight or loose?
We believe it is loose because
$\Delta(4k-1) + (2k+1) \Delta(3) = \Delta(3k+1)+2$,
and a class of the form ``$\Delta(.)+const$''
is always loose (we think).
Incidentally,
the number of the triangles in the class is
unbounded with $k$.

However,
if (as we believe)
the sequence of densest packings $n_p (k), k=1,2,3,...$ for 
fixed $p=2$ can be continued with all packings having the same pattern
of six triangles,
then t58a might also be a member of the class 
$\Delta(2k+1) + 5\Delta(k)$ for $k=3$.
The next term in the latter class would be the
densest packing of $n=n_2 (4)=\Delta(9) + 5\Delta(4)=95$ disks.
Our experiments with 95 disks, indeed, produced
the desired pattern of six triangles in the densest packing t95a
(Fig.~\ref{t58-175}).
This reinforces our suspicion that the class $n_2 (k),~k=1,2,...$, exists
in which each densest
packing consists of six densely packed triangles.
The class $n_2 (k),~k=1,2,...$, should be tight
because each densest packing in it consists of a fixed number of triangles.

By increasing $p$ we are moving into a different class,
which is again tight if the conjecture above is correct.
Thus, for $p=3$ we have the sequence 
$n_3 (1)=22, n_3 (2) = 57, n_3 (3) = 108, n_3 (4) = 175,....$
Packings t22a (Fig~\ref{t22}), t57a (Fig~\ref{t57-93}) indeed
each consist of $2(p+1)=8$ densely packed triangles.
Our experiments with 108 and 175 disks yield
the same pattern in the best packings (Fig~\ref{t58-175})
so the class $n_3 (k)$ probably exists too.
In the same way,
the class $n_p (k), k=1,2,...$ exists for any fixed index
$p$ and has a distinct pattern with $r=2(p+1)$ triangles and $p-1$ rattlers.

If this is correct,
then Figure \ref{t58-175} can be seen as the $2 \times 2$ submatrix
for $2 \le p \le 3$ and $3 \le k \le 4$
of the matrix of dense packings of $n_p (k)$ disks
where $1 \le k,p \le \infty$.
By traversing a row or a column of this matrix we obtain a distinct
infinite class of packings.
Our conjecture is that each row class is tight and each column class is loose.

This matrix contains three infinite classes conjectured in Section~5:
the row at $p = 1$ is the class $4 \Delta(k)$,
the column at $k = 1$ is the class $\Delta (2p) +1$,
and the column at $k = 2$ is the class 
$\Delta (3p+1) +2 = \Delta(3p-1) + (2p+1) \Delta(2)$.

The conjecture about the full matrix
is also reinforced by the fact
that $k$ and $p$ in (3) are unique for each given value of $n = n_p (k)$.
This can be easily seen using the first expression for $n_p (k)$
in (3).
To further test our matrix conjecture,
we generated the list of all $n$ of the form $n=n_p (k)$
for $n \le 300$.
These are:
4,  11, 12, 22, 24, 30, 37, 40, 56, 57, 58, 60, 79, 84, 93, 95, 106, 108, 
112, 137, 138, 141, 144, 172, 174, 175, 180, 192, 196, 211, 220, 254, 255, 
256, 258, 260, 264, and 280.
Some increments in this increasing sequence are small.
Specifically, in each following subsequence increments do not exceed 2:
(11, 12), (22, 24), (56, 57, 58, 60), (93, 95), (106, 108), (137, 138),
(172, 174, 175), and (254, 255, 256, 258, 260).

Now, take for example,
$255 = n_7 (2)$ and $256 = n_5 (3)$.
These are two ``almost'' equal numbers of disks.
However according to the matrix conjecture they should produce
different patterns of densest packings:
the pattern for 255 should consist of one large and 15 small triangles
with 6 rattlers and the pattern for 256 of one large and 11 small
triangles and 4 rattlers.
Similarly, the matrix conjecture prescribes specific patterns
for the densest packing of the other numbers of disks $n$ 
of this sequence, e.g., $254 = n_{11} (1)$,
$258 = n_3 (5)$, and $260 = n_2 (7)$.

One might think it would be a stress test for both the matrix conjecture
and our packing procedure to try to pack these numbers of disks.
Note that many of the packings for smaller values
in the sequences above have been generated
(as discussed above)
and they all conform to the matrix conjecture.
Thus, we experimented with packing
$n=$ 254, 255, 256, 258, and 260 disks.
Recall that the procedure of packing
has no idea, so to say, of the desirable packing.
Starting with random initial conditions the disks perform
chaotic movements, they collide with each other
and with the boundaries millions of times
and 
each collision evaluation is subject
to roundoff error.

\begin{figure}[H]
\centerline{\psfig{file=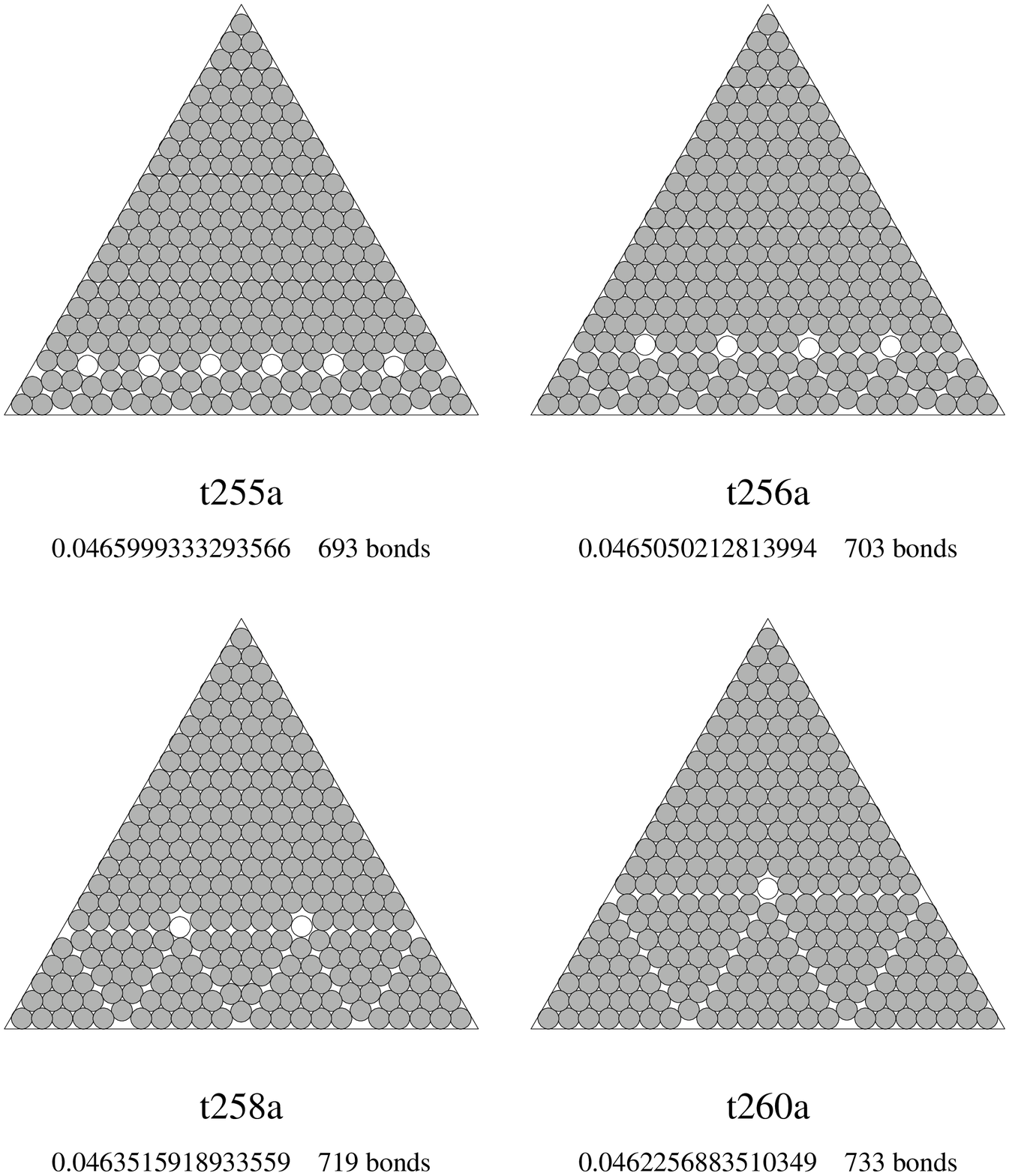,width=7.1in}}

\caption{The best packings of 255 disks (t255a), 256 disks (t256a), 258 disks (t258a), and 260 disks (t260a).}
\label{t255-60}
\end{figure}

The experiments turned out to be not so difficult.
(Case of 53 disks proved to be harder.)
As expected, the best packing of 254 disks
has the pattern of its class $\Delta (2k) +1$ with 10 rattlers
($d(254)=0.0467170396481042$, 679 bonds; we omit the diagram).
Fig.\ref{t255-60} shows the patterns of
packings t255a, t256a, t258, and t260a.
These too are consistent with formula (3).
\addtolength{\textheight}{+.25in}
Note that because of the large number of disks in the packings
the scale of drawing in Fig.\ref{t255-60} is small
and bonds are not seen.
The pictures with a larger scale (omitted here) show 
that all the bonds exist in right places.

\clearpage
\section{Discussion}
\hspace*{\parindent}
While a finite number of patterns for infinite classes have been 
tentatively identified to date 
(two one-parameter patterns known or conjectured previously
joined by several such patterns in Sec.5
and a two-parameter ``matrix'' pattern in Sec.6)
a countable infinity of such patterns and classes probably exists.
Furthermore, each value of $n$ may well be a member of one or more
such classes.
Thus, the values $n=18$, 23, 32, and 33, which
were not placed into classes in this paper,
may well be members of as yet unidentified classes of packings 
with complex patterns.  In fact, a fixed value of $n$
may be on the paths of many, possibly infinitely many, such classes.
12 disks gives an example of this: it is on the path
of the class $4\Delta(k)$ and it is also the first term
of the classes $\Delta(3k+1)+2$ and $\Delta(2k+3)-3$.
As the value of $n$ increases along the path of a class, 
``hesitations'' of the best pattern may occur,
wherein several different nonequivalent patterns
coexist among the rigid packings and compete 
for the title of the best.  A resolved case of such hesitation
occurs for the class $\Delta(2k+1)+1$ where for $k \ge 5$ ($n \ge 67$)
two equivalent best patterns finally emerge 
(at least according to our experiments).
We were not able to confirm by experiments the winning
pattern for the class $\Delta(k+2)-2$.
Will such hesitation always be resolved in favor of one
of the competing patterns in a finite initial
segment of the path?


\begin{thebibliography}{MMMM}
\bibitem[CFG]{CFG}
H. T. Croft, K. J. Falconer and R. K. Guy,
{\em Unsolved Problems in Geometry}, Springer Verlag, Berlin, 1991, 107--111.
\bibitem[FG]{FG}
J.~H. Folkman and R.~L. Graham,
A packing inequality for compact convex subsets of the plane,
{\em Canad. Math. Bull.} {\bf 12} (1969), 745--752.
\bibitem[L]{L}
B. D. Lubachevsky,
How to simulate billiards and similar systems, {\em J. Computational
Physics} {\bf 94} (1991), 255--283.
\bibitem[LS]{LS}
B. D. Lubachevsky and F. H. Stillinger,
Geometric properties of random disk packings,
{\em J. Statistical Physics} {\bf 60} (1990), 561--583.
\bibitem[M1]{M1}
J.~B.~M. Melissen,
Densest packings for congruent circles in an equilateral triangle,
{\em Amer. Math. Monthly} {\bf 100} (1993), 916--925.
\bibitem[M2]{M2}
J.~B.~M. Melissen, 
Optimal packings of eleven equal circles in an equilateral triangle,
{\em Acta Math. Hung.} {\bf 65} (1994), 389--393.
\bibitem[MS]{MS}
J.~B.~M. Melissen and P.~C. Schuur,
Packing 16, 17 or 18 circles on an equilateral triangle,
{\em Disc. Math.} (to appear).
\bibitem[N]{N}
D. J. Newman, private communication.
\bibitem[O]{O}
N. Oler,
A finite packing problem, {\em Canad. Math. Bull.} {\bf 4} (1961), 153--155.
\end{thebibliography}
\end{document}